\input amstex
\documentstyle{amsppt}
\magnification=\magstep1 \NoRunningHeads

\topmatter

\title
Mixing actions of Heisenberg group
 \endtitle
\author  Alexandre~I.~Danilenko
\endauthor
\abstract
Mixing (of all orders) rank-one actions $T$ of Heisenberg group $H_3(\Bbb R)$ are constructed.
The restriction of $T$ to the center of $H_3(\Bbb R)$
is simple and commutes only with $T$.
Mixing Poisson and mixing Gaussian actions of $H_3(\Bbb R)$ are also constructed.
A rigid weakly mixing rank-one action $T$ is constructed such that
the restriction of $T$ to the center of $H_3(\Bbb R)$ is not isomorphic to its inverse.
\endabstract

\address
 Institute for Low Temperature Physics
\& Engineering of National Academy of Sciences of Ukraine, 47 Lenin Ave.,
 Kharkov, 61164, UKRAINE
\endaddress
\email alexandre.danilenko\@gmail.com
\endemail

\NoBlackBoxes
\endtopmatter
\document

\head 0. Introduction
\endhead

We  state a question by Dan Rudolph:

\proclaim{Problem 1} Which amenable locally compact second countable groups $G$ admit mixing free actions with zero entropy?
\endproclaim

Consider a subtler problem:

\proclaim{Problem 2} Which amenable locally compact second countable groups $G$ admit mixing rank-one free actions?
\endproclaim

We recall that a measure preserving action $T=(T_g)_{g\in G}$
of $G$ on a standard probability space $(X,\goth B,\mu)$ is
\roster
\item"---" {\it mixing} if
$\lim_{g\to G}\mu(T_gA\cap B)=\mu(A)\mu(B)$ for all $A,B\in\goth B$,
\item"---" {\it mixing of order $l$} if for each $\epsilon>0$ and subsets
$A_0,A_1,\dots,A_l\in\goth B$, there is a compact subset $K\subset G$ such that
$$
|\mu(T_{g_0}A_0\cap\dots\cap T_{g_l}A_l)-\mu(A_0)\mu(A_1) \cdots\mu(A_l)|<\epsilon
$$
for each collection $g_0,\dots g_l\in G$ with $g_ig_j^{-1}\not\in K$,
\item"---" {\it rigid} if there is a sequence $g_n\to\infty$ in $G$ such that
$\lim_{g\to G}\mu(T_gA\cap B)=\mu(A\cap B)$ for all $A,B\in\goth B$
\endroster

We recall now the definition of  rank one.
Fix a F{\o}lner sequence $(F_n)_{n=1}^\infty$ in $G$.

\roster
\item"---" A {\it Rokhlin tower or column} for $T$ is a triple $(Y,f, F)$,
 where $Y\in\frak B$, $F$ is a relatively compact subset of $G$
and $f:Y\to F$ is a measurable mapping
 such that for any Borel subset $H\subset F$ and an element $g\in G$ with
 $gH\subset F$, one has
 $f^{-1}(gH)=T_gf^{-1}(H)$.
\item"---"
The action  $T$ has {\it rank one along $(F_n)_{n\in\Bbb N}$} if
there exists a sequence of Rokhlin towers $(Y_n,f_n,F_n)$ such that
  for each subset $B\in\frak B$, there is a sequence
of Borel subsets $H_n\subset F_n$ such that
$$
\lim_{n\to\infty}\mu(B\triangle f_n^{-1}(H_n))=0.
$$
\endroster

It is easy to see that each  rank-one action is ergodic.
It is well known that each rank-one action has  zero entropy.
Hence every solution of Problem~2 is a solution of Problem~1.
The two problems are still open.
However  various constructions of mixing rank-one actions are elaborated for the following amenable groups:
\roster
\item"---" $G=\Bbb Z$ in \cite{Or}, \cite{Ad}, \cite{CrSi}, \cite{Ry4}, etc.
\item"---" $G=\Bbb Z^2$ in \cite{AdSi},
\item"---"
$G=\Bbb R$ in \cite{Pr}, \cite{Fa},
\item"---"
$G=\Bbb R^{d_1}\times\Bbb Z^{d_2}$ for arbitrary $d_1,d_2\ge 0$ in \cite{DaSi1},
\item"---"
$G=\bigoplus_{j=0}^\infty G_j$, where $G_j$ is a finite group \cite{Da2},
\item"---"
$G$
is a locally normal discrete countable group, i.e. $G=\bigcup_{j=1}^\infty G_j$ for a nested sequence $G_1\subset G_2\subset\cdots$ of normal finite subgroups $G_j\subset G$ \cite{Da4}.
\endroster

The Heisenberg group $H_3(\Bbb R)$ consists of of $3\times 3$ upper triangular matrices of the form
$$
    \pmatrix
 1 & a & c\\ 0 & 1 & b\\ 0 & 0 & 1\\
\endpmatrix,
$$
where $a,b,c$ are arbitrary reals.
The Heisenberg group endowed with the natural topology is a connected, simply-connected nilpotent Lie group.
We now let
$$
a(t):=
\pmatrix
 1 & t & 0\\ 0 & 1 & 0\\ 0 & 0 & 1\\
\endpmatrix,
\quad
b(t):=
\pmatrix
 1 & 0 & 0\\ 0 & 1 & t\\ 0 & 0 & 1\\
\endpmatrix,
\quad
c(t):=\pmatrix
 1 & 0 & t\\ 0 & 1 & 0\\ 0 & 0 & 1\\
\endpmatrix.
$$
Then $\{a(t)\mid t\in\Bbb R\}$, $\{b(t)\mid t\in\Bbb R\}$ and $\{ c(t)\mid t\in\Bbb R\}$ are three closed one-parameter subgroups in $H_3(\Bbb R)$.
The latter one is the center of $H_3(\Bbb R)$.
Every element  $g$ of $H_3(\Bbb R)$  can be written uniquely as the product $g=a(t_1)b(t_2)c(t_3)$ for some $t_1,t_2,t_3\in\Bbb R$.
We also note that $[a(t_1),b(t_2)]:=a(t_1)b(t_2)a(t_1)^{-1}b(t_2)^{-1}=c(t_1t_2)$.

We now state one of the main results of the present paper.

\proclaim{Theorem~0.1} There exist mixing of all orders rank-one (and hence zero entropy) action
$T$
of Heisenberg group $H_3(\Bbb R)$.
\endproclaim

We construct these actions utilizing the $(C,F)$-construction with randomly chosen `spacers' in the spirit of Ornstein's rank-one mixing map \cite{Or} (see also \cite{Ru}).
This construction is an algebraic counterpart of the well-known  inductive construction process of `cutting-and-stacking with a single tower' for rank-one maps.
It was introduced in \cite{dJ} (see also \cite{Da1} and a survey \cite{Da3}) as a convenient tool to produce rank-one actions of groups with torsions or non-Abelian groups.
We show first that the restriction of $T$ to the center of $H_3(\Bbb R)$ is mixing.
Then we adapt  Ryzhikov's idea from \cite{Ry1} to deduce that the entire action is mixing too.

Since the discrete countable Heisenberg group $H_3(\Bbb Z)$ is a co-compact lattice in $H_3(\Bbb Z)$\footnote{$H_3(\Bbb Z)$ is defined in a similar way as $H_3(\Bbb R)$ but with $a,b,c\in\Bbb Z$.}, the restriction of a mixing zero-entropy action  of $H_3(\Bbb R)$ to $H_3(\Bbb Z)$ yields a mixing  zero-entropy action of $H_3(\Bbb Z)$.
Thus it follows from Theorem~0.1 that $H_3(\Bbb R)$ and $H_3(\Bbb Z)$ are both in the list of solutions of Problem~1.

In order to state the next problem considered in this paper we recall some definitions from the theory of joining of dynamical systems (see \cite{dJRu}, \cite{Gl}, \cite{Th}, \cite{dR}).
Given an ergodic action $T=(T_g)_{g\in G}$ of a  locally compact second countable group $G$ on a standard probability space $(X,\goth B,\mu)$, we denote by $C(T)$ the centralizer of $T$, i.e. the set of $\mu$-preserving invertible transformations of $X$ commuting with $T_g$ for each $g\in G$.
Given a transformation $R\in C(T)$, we denote by $\mu_R$ the corresponding off-diagonal measure on $X\times X$, i.e.
$\mu_R(A\times B):=\mu(RA\cap B)$.
A $(T_g\times T_g)_{g\in G}$-invariant measure on $X\times X$ whose  marginal on each copy of $X$ is $\mu$ is called a {\it 2-fold self-joining} of $T$.
For $l>1$,  $l$-fold self-joinings of $T$ are defined in a similar way.
Let $J_2^e(T)$ denote the set of ergodic 2-fold self-joinings of $T$.
If $J_2^e(T)\subset\{\mu_R\mid R\in C(R)\}\cup\{\mu\times\mu\}$ then
$T$ is called 2-fold simple.
If $T$ is 2-fold simple and for each $l>1$, the $l$-fold Cartesian power of $\mu$  is the only  $l$-fold self-joining of $T$ whose projection on the product of any two copies of $X$ in $X^l$ equals $\mu\times \mu$ then
$T$ is called {\it simple}.
If $T$ is simple and $C(T)\subset\{T_g\mid g\in G\}$ then $T$ is said to have the  property of {\it minimal self-joinings} (MSJ).

In \cite{dJ}, A. del Junco raised  the following problem.

\proclaim{Problem 3}
Given a locally compact second countable group $G$ and a closed non-compact subgroup $L\subset G$, is there a free action $T$ of $G$
such that the sub-action $(T_l)_{l\in L}$ is weakly mixing and 2-fold simple and the centralizer of this sub-action is $\{T_h\mid h\in C_G(L)\}$, where $C_G(L)$ stands for the centralizer of $L$ in $G$,
i.e. $C_G(L)=\{g\in G\mid gl=lg\text{ for all }l\in L\}$.
\endproclaim

 He shows that for some special pairs $L\subset G$,  solutions of Problem~3 lead to non-trivial counterexamples in ergodic theory \cite{dJ}.
In the most interesting cases  $L$ is the center of $G$ and hence $C_G(L)=G$.
We now write the second main result of the paper.

\proclaim{Theorem 0.2} There exists rank-one mixing action $T$ of $H_3(\Bbb R)$  such that the following are satisfied.
\roster
\item"\rom{(i)}"
The flow $(T_{c(t)})_{t\in\Bbb R}$ is simple
and $C((T_{c(t)})_{t\in\Bbb R})=\{T_g\mid g\in H_3(\Bbb R)\}$.
\item"\rom{(ii)}"
The transformation $T_{c(1)}$ is simple
and $C(T_{c(1)})=\{T_g\mid g\in H_3(\Bbb R)\}$.
\item"\rom{(ii)}" $T$ has MSJ.
\endroster
\endproclaim

Thus Theorem~0.2 answers Problem~3 affirmatively in the particular case when $G=H_3(\Bbb R)$ and $L$ is either the center of $H_3(\Bbb R)$ or $L=\{c(n)\mid n\in\Bbb Z\}$.

Next, we construct mixing rank-one infinite measure preserving actions of Heisenberg group within the framework of $(C,F)$-construction.
We recall that an infinite measure preserving action $T$ of $H_3(\Bbb R)$ is called {\it mixing} (or {\it $0$-type}) if $\mu(T_gA\cap B)\to 0$ as $g\to\infty$ for all subsets $A,B\subset X$ of finite measure.
Mixing actions in infinite measure need not be ergodic (see \cite{DaSi2} and references therein) however the rank one  implies ergodicity.
We note that while the construction of the finite measure preserving mixing actions of $H_3(\Bbb R)$ in Theorem~0.1 is of stochastic  nature (as in \cite{Or}, \cite{Ru}, \cite{dJ}),
all the parameters of the infinite measure preserving mixing actions of $H_3(\Bbb R)$ are determined {\it explicitly} (effectively).
 For that we apply a technique of {\it fast-growing spacers} suggested by V.~Ryzhikov in \cite{Ry4} for  the case of $\Bbb Z$-actions.
As an application we obtain the following corollary.

\proclaim{Corollary 0.3} There exist mixing Poisson and mixing Gaussian (probability preserving) actions of $H_3(\Bbb R)$.
\endproclaim

We also consider a problem of {\it asymmetry} for ergodic actions of Heisenberg group.

\proclaim{Theorem 0.4}
There is a rigid weakly mixing rank-one action $T$ of $H_3(\Bbb R)$
such that the transformation $T_{c(1)}$ is ergodic and non-conjugate to its inverse $T_{c(1)}^{-1}$.
\endproclaim

To prove this theorem we `incorporate' V.~Ryzhikov's example of an asymmetric rank-one transformation from \cite{Ry3} (see also \cite{DaRy} for a similar example of an asymmetric rank-one flow) into the construction of mixing actions of $H(\Bbb R)$ from~Theorem~0.1.

The outline of the paper is the following.
Sections~1 and 2 are  preliminary.
We describe the unitary dual of $H_3(\Bbb R)$ and state the spectral decomposition theorem for unitary representations of $H_3(\Bbb R)$ in Section~1.
In Section~2 we remind the $(C,F)$-construction for locally compact group actions.
Sections~3 and 4 are devoted to the proof of Theorems~0.1 and 0.2 respectively.
In Section~5 we construct a
mixing rank-one infinite measure preserving action of $H_3(\Bbb R)$
and deduce Corollary~0.3.
Section~6 is devoted to asymmetric actions of $H_3(\Bbb R)$.
We prove there Theorem~0.4.
In Section~7 we apply the spectral decomposition to the Koopman representations of $H_3(\Bbb R)$.
Section~8 consists of concluding remarks and open problems.

\head 1. Heisenberg group and its unitary dual
\endhead

The Lie algebra of $H_3(\Bbb R)$ is
$$
\goth h_3:=\left\{\pmatrix
 0 & a & c\\ 0 & 0 & b\\ 0 & 0 & 0\\
\endpmatrix \Bigg| \,a,b,c\in\Bbb R\right\}.
$$
We endow it with the natural  topology.
Then the exponential map $\exp:\goth h_3\to H_3(\Bbb R)$ is a homeomorphism.

The subgroups $H_{2,a}:=\{a(t_1)c(t_3)\mid t_1,t_3\in\Bbb R\}$
and $H_{2,b}:=\{b(t_2)c(t_3)\mid t_2,t_3\in\Bbb R\}$ are both Abelian,  normal and closed in $H_3(\Bbb R)$.
The corresponding group extensions
$$
\align
&0\to H_{2,a}\to H_3(\Bbb R)\to H_3(\Bbb R)/H_{2, a}\to 0\quad \text{and}\\
&0\to H_{2,b}\to H_3(\Bbb R)\to H_3(\Bbb R)/H_{2, b}\to 0
\endalign
$$
are both split.
 This implies that $H_3(\Bbb R)$ is isomorphic to  the semidirect product $\Bbb R^2\rtimes_A\Bbb R$, where the homomorphism $A:\Bbb R\to \text{GL}_2(\Bbb R)$ is given by $A(t):=\pmatrix 1&t\\
0&1\endpmatrix$, $t\in\Bbb R$.
The subgroups $H_{2,a}$ and $H_{2,b}$ are automorphic in $H_3(\Bbb R)$,
i.e. there is an isomorphism $\theta$ of $H_3(\Bbb R)$ with $\theta(H_{2,a})=H_{2,b}$.
We define $\theta$ by setting $\theta(a(t)):=b(t)$, $\phi(b(t)):=a(t)$ and $\theta(c(t)):=c(-t)$ for all $t\in\Bbb R$.
To put it in other way,
$$
\theta
\pmatrix
 1 & a & c\\ 0 & 1 & b\\ 0 & 0 & 1\\
\endpmatrix=
\pmatrix
 1 & b & ab-c\\ 0 & 1 & a\\ 0 & 0 & 1\\
\endpmatrix.
$$
We call $\theta$ the {\it flip} in  $H_3(\Bbb R)$.
We note that $\theta^2=\text{id}$.

The set of unitarily equivalent classes of irreducible (weakly continuous) representations of
$H_3(\Bbb R)$ is called the {\it unitary dual} of
$H_3(\Bbb R)$.
It is denoted by $\widehat {H_3(\Bbb R)}$.
The irreducible unitary representations of $H_3(\Bbb R)$ are well known.
They consist of  a family of 1-dimensional representations $\pi_{\alpha,\beta}$, $\alpha,\beta\in\Bbb R$ and a family of infinite dimensional representations $\pi_\gamma$, $\gamma\in\Bbb R\setminus\{0\}$ as follows \cite{Ki}:
$$
\align
\pi_{\alpha,\beta}(c(t_3)b(t_2)a(t_1))&:=e^{i(\alpha t_1+\beta t_2)} \quad\text{and }\\
(\pi_\gamma(c(t_3)b(t_2)a(t_1))f)(x)&:=e^{i\gamma (t_3+t_2x)}f(x+t_1),\quad f\in L^2(\Bbb R,\lambda_{\Bbb R}).
\endalign
$$
Thus we can identify $\widehat{H_3(\Bbb R)}$  with the disjoint union $\Bbb R^2\sqcup \Bbb R^*$.
There is a natural Borel $\sigma$-algebra on the  unitary dual of each locally compact second countable group \cite{Ma}.
 In the case of the Heisenberg group this Borel $\sigma$-algebra coincides with the standard $\sigma$-algebra of Borel subsets in  $\Bbb R^2\sqcup\Bbb R^*$.
Given an arbitrary unitary representation $U=(U(g))_{g\in H_3(\Bbb R)}$
of $H_3(\Bbb R)$ in a separable Hilbert space $\Cal H$, there are  a measure $\sigma_U$ on $\widehat{H_3(\Bbb R)}$ (i.e. two measures $\sigma_U^{1,2}$ on $\Bbb R^2$ and $\sigma_U^3$ on $\Bbb R^*$) and a map $l_U:\widehat{H_3(\Bbb R)}\to\Bbb N\cup\{\infty\}$  (i.e. two maps $l_U^{1,2}:\Bbb R^2\ni (x,y)\mapsto l_U^{1,2}(x,y)\in\Bbb N\cup\{\infty\}$
 and
$l_U^3:\Bbb R^*\ni z\mapsto l_U^3(z)\in\Bbb N\cup\{\infty\}$) such that the following decompositions hold up to unitary equivalence
$$
\align
\Cal H&= \int^{\oplus}_{\Bbb R^2}\,\bigoplus_{j=1}^{l_U^{1,2}(\alpha,\beta)}\Bbb C\,d\sigma_U^{1,2}(\alpha,\beta)
\oplus\int_{\Bbb R^*}^\oplus\bigoplus_{j=1}^{l_U^3(\gamma)}L^2(\Bbb R,\lambda_\Bbb R)\,d\sigma_U^3(\gamma)\quad\text{and} \\
U(g)&=\int^{\oplus}_{\Bbb R^2}\,\bigoplus_{j=1}^{l_U^{1,2}(\alpha,\beta)} \pi_{\alpha,\beta}(g)\,d\sigma_U^{1,2}(\alpha,\beta)
\oplus\int_{\Bbb R^*}^\oplus\bigoplus_{j=1}^{l_U^3(\gamma)}\pi_\gamma(g)\,d\sigma_U^3(\gamma).
\endalign
$$
The equivalence class of  $\sigma_U$ is called the {\it maximal spectral type} of $U$.
The map $l_U$ is called the  {\it multiplicity function} of $U$.
The essential range if $l_U$ is called the {\it set of spectral multiplicities} of $U$.
The maximal spectral type
and the multiplicity function of $U$ ($\sigma_U$-mod 0)
are both determined uniquely by the unitarily equivalent class of $U$.

Given a measure preserving action $T$ of $H_3(\Bbb R)$ on a standard probability space $(X,\goth B,\mu)$, we can associate  a unitary representation $U_T$ of of $H_3(\Bbb R)$ in $L^2(X,\mu)$ by setting
$U_T(g)f:=f\circ T_g^{-1}$.
It is called the Koopman representation of $H_3(\Bbb R)$ associated with $T$.
The  maximal spectral type of $U_T$ is called the maximal spectral type of $T$ and the maximal spectral type of the restriction of $U_T$ to the subspace of zero mean functions in $L^2(X,\mu)$ is called the restricted maximal spectral type of~$T$.

\head 2. $(C,F)$-construction
\endhead

Let $G$ be a unimodular l.c.s.c. amenable group.
Fix a $\sigma$-finite  Haar measure $\lambda_G$ on it.
Given two subsets $E,F\subset G$, by
$EF$ we mean their algebraic product, i.e. $EF=\{ef\mid e\in E, f\in F\}$. The set
$\{e^{-1}\mid e\in E\}$ is denoted by $E^{-1}$.
If $E$ is a singleton, say
$E=\{e\}$, then we will write $eF$ for $EF$ and $Fe$ for $FE$.

To define a $(C,F)$-action of $G$ we need two sequences $(F_n)_{n\ge 0}$
and $(C_n)_{n>0}$ of subsets in $G$ such that the following conditions are
satisfied:

\roster
\item"(I)"
$(F_n)_{n=0}^\infty\text{ is a F{\o}lner sequence in }G$,
\item"(II)"
$C_n\text{ is finite and }\# C_n>1$,
\item"(III)"
$F_nC_{n+1}\subset F_{n+1}$,
\item"(IV)"
 $F_nc\cap F_nc'=\emptyset\text{ for all }c\ne c'\in C_{n+1}$.
\endroster
 We
put $X_n:=F_n\times\prod_{k>n}C_k$, endow $X_n$ with the standard Borel
product $\sigma$-algebra and define a Borel embedding $X_n\to X_{n+1}$ by
setting
$$
(f_n,c_{n+1},c_{n+2},\dots)\mapsto (f_nc_{n+1}, c_{n+2},\dots).\tag2-1
$$
It is well defined due to (III) and (IV).
Then we have $X_1\subset
X_2\subset\cdots$.
Hence $X:=\bigcup_nX_n$ endowed with the natural Borel
$\sigma$-algebra, say $\frak B$, is a standard Borel space. Given a
Borel subset $A\subset F_n$, we put
$$
[A]_n:=\{x\in X\mid x=(f_n,c_{n+1}, c_{n+2}\dots)\in X_n\text{ and } f_n\in
A\}
$$
and call this set an $n$-{\it cylinder}.
It is clear that the
$\sigma$-algebra $\frak B$ is generated by the family of all cylinders.

Now we are going to define a `canonical' measure on $(X,\frak B)$.
Let
$\kappa_n$ stand for the equidistribution on $C_n$.
Let
 $\nu_n:= (\#C_1
\cdots \#C_n)^{-1}\lambda_G \upharpoonright
 F_n$.
We define an infinite  product
measure $\mu_n$ on $X_n$ by setting
$
\mu_n:=\nu_n\times\kappa_{n+1}\times\kappa_{n+2}\times\cdots,
$
$n\in\Bbb N$.
Then the embeddings \thetag{2-1} are all measure preserving.

Hence a $\sigma$-finite measure $\mu$ on $X$ is well defined by the
restrictions $\mu\restriction X_n:=\mu_n$, $n\in\Bbb N$.
To put it in
another way, $(X,\mu)=\injlim_n(X_n,\mu_n)$.
Since
$$
\mu_{n+1}(X_{n+1})=\frac{\nu_{n+1}(F_{n+1})}{\nu_{n+1}(F_nC_{n+1})}
\mu_n(X_n)=\frac{\lambda_G(F_{n+1})}{\lambda_G(F_n)\# C_{n+1}}\mu_n(X_n),
$$
it follows that $\mu$
 is finite if and only if
$$
 \prod_{n=0}^{\infty}\frac{\lambda_G(F_{n+1})}{\lambda_G(F_n)\#
C_{n+1}} <\infty,\text{ i.e.,
}\sum_{n=0}^\infty\frac{\lambda_G(F_{n+1}\setminus(F_nC_{n+1}))}
{\lambda_G(F_n)\# C_{n+1}}<\infty.
\tag2-2
$$
To construct a $\mu$-preserving action of $G$ on $(X,\mu)$, we fix a
filtration $K_1\subset K_2\subset\cdots$ of $G$ by compact subsets. Thus
$\bigcup_{m=1}^\infty K_m=G$. Given $n,m\in \Bbb N$, we set
$$
\align
 D_m^{(n)}&:=\bigg(\bigcap_{k\in K_m}(k^{-1}F_n)\cap
F_n\bigg)\times\prod_{k>n}C_k\subset X_n\text{ and
}\\
R_m^{(n)}&:=\bigg(\bigcap_{k\in K_m}(kF_n)\cap
F_n\bigg)\times\prod_{k>n}C_k\subset X_n.
\endalign
$$
It is easy to verify that
$
D_{m+1}^{(n)}\subset D_m^{(n)}\subset D_m^{(n+1)}
 \text{ and }R_{m+1}^{(n)}\subset R_m^{(n)}\subset R_m^{(n+1)}.
$
We define a Borel mapping
$
K_m\times D_m^{(n)}\ni(g,x)\mapsto T_{m,g}^{(n)}x\in R_m^{(n)}
$
by setting for $x=(f_n,c_{n+1},c_{n+2},\dots)$,
$$
T^{(n)}_{m,g}(f_n,c_{n+1},c_{n+2}\dots):=(gf_n,c_{n+1}, c_{n+2},\dots).
$$
Now let $D_m:=\bigcup_{n=1}^\infty D_m^{(n)}$ and
$R_m:=\bigcup_{n=1}^\infty R_m^{(n)}$. Then a Borel mapping
$$
T_{m,g}:K_m\times D_m\ni(g,x)\mapsto T_{m,g}x\in R_m
$$
is well defined by the restrictions $T_{m,g}\restriction
D_m^{(n)}:=T_{m,g}^{(n)}$ for $g\in K_m$ and $n\ge 1$.
It is easy to see
that $D_m\supset D_{m+1}$, $R_m\supset R_{m+1}$ and $T_{m,g}\restriction
D_{m+1}=T_{m+1,g}$ for all $m$.
It follows from (I) that
$\mu_n(X_n\setminus D_m^{(n)})\to 0$ and $\mu_n(X_n\setminus R_m^{(n)})\to 0$ as $n\to\infty$. 
Hence
$\mu(X\setminus D_m)=\mu(X\setminus R_m)=0$ for all $m\in\Bbb N$.
Finally we set $\widehat
X:=\bigcap_{m=1}^\infty D_m\cap\bigcap_{m=1}^\infty R_m$ and define a Borel
mapping $ T:G\times \widehat X \ni(g,x)\to T_gx\in \widehat X $ by setting
$T_gx:=T_{m,g}x$ for some (and hence any) $m$ such that $g\in K_m$.
It is
clear that $\mu(X\setminus \widehat X)=0$.
Thus, we obtain that $T=(T_g)_{g\in G}$ is
a free Borel measure preserving action of $G$ on a $\mu$-conull subset of the
standard $\sigma$-finite space $(X,\frak B,\mu)$. 
It is easy to see that
$T$ does not depend on the choice of filtration $(K_m)_{m=1}^\infty$.
Throughout the paper we do not distinguish between two measurable sets (or
mappings) which agree almost everywhere.

\definition{Definition 2.1}
$T$ is called the $(C,F)$-{\it action of $G$
 associated with} $(C_{n+1},F_n)_{n=0}^{\infty}$.
 \enddefinition

We now list some basic properties of $(X,\mu,T)$.
Given Borel subsets
$A,B\subset F_n$, we have
$$
\align
&[A\cap B]_n=[A]_n\cap[B]_n,\ [A\cup B]_n=[A]_n\cup[B]_n,\tag2-3\\
&[A]_n=[AC_{n+1}]_{n+1}=\bigsqcup_{c\in C_{n+1}}[Ac]_{n+1},\tag2-4\\
&T_g[A]_n=[gA]_n\text{ if } gA\subset F_n,\tag2-5\\
&\mu([A]_n)=\# C_{n+1}\cdot\mu([Ac]_{n+1})\text{ for every }c\in
C_{n+1},\tag2-6\\
&\mu([A]_n)=\frac{\lambda_G(A)}{\lambda_G(F_n)}\mu(X_n).\tag2-7
\endalign
$$

 In case
$G=\Bbb Z$, it is easy to notice a similarity between the
$(C,F)$-construction and the classical cutting-and-stacking construction of
rank-one transformations. 
Indeed, $F_{n-1}$ (or, more
precisely, the set of $(n-1)$-cylinders) corresponds to the levels of the
$(n-1)$-tower and $C_n$ corresponds to the locations of the copies of
$F_{n-1}$ inside the $n$-th tower $F_n$. (The copies $F_{n-1}c$, $c\in
C_n$, are disjoint by (IV) and they sit inside $F_n$ by (III).)
The remaining part of $F_n$, i.e. $F_n\setminus(F_{n-1}C_n)$, is the set of
spacers in the $n$-th tower.

Each
$(C,F)$-action is of  rank one.

\head 3. Mixing rank-one actions of Heisenberg group
\endhead

Given three positive reals $\alpha,\beta$ and $\gamma$, we set
$$
\align
I(\alpha,\beta,\gamma)
&:=\{c(t_3)b(t_2)a(t_1)\mid |t_1|\le\alpha,
|t_2|\le\beta,|t_3|\le\gamma\}\\
&=\left\{
\pmatrix
1 & t_1 & t_3\\ 0 & 1 & t_2\\ 0 & 0 & 1\\
\endpmatrix
\bigg| \, |t_1|\le\alpha,
|t_2|\le\beta,|t_3|\le\gamma
\right\}.
\endalign
$$
It is easy to verify that
$$
\align
I(\alpha,\beta,\gamma)I(\alpha',\beta',\gamma') &\subset I(\alpha+\alpha',\beta+\beta',\gamma+\gamma'+\alpha\beta')\quad \text{and}\\
I(\alpha,\beta,\gamma)^{-1} &\subset I(\alpha,\beta,\gamma+\alpha\beta).
\endalign
$$
Now we define a map $\phi_{\alpha,\beta,\gamma}:\Bbb Z^3\to H_3(\Bbb R)$ by setting
$$
\phi_{\alpha,\beta,\gamma}(j_1,j_2,j_3) :=c(2\gamma j_3)b(2\beta j_2)a(2\alpha j_1), \quad j_1,j_2,j_3\in\Bbb Z.
$$
This map is not a group homomorphism.
However for all $j_1,j_2,j_3,p\in\Bbb Z$, we have
$$
\align
\phi_{\alpha,\beta,\gamma}(0,0,p)\phi_{\alpha,\beta,\gamma}(j_1,j_2,j_3)
&=
\phi_{\alpha,\beta,\gamma}(0,0,p)\phi_{\alpha,\beta,\gamma}(j_1,j_2,j_3)\\
&=
\phi_{\alpha,\beta,\gamma}(j_1,j_2,p+j_3).
\endalign
$$
The following tiling property holds:
$$
H_3(\Bbb R)=\bigsqcup_{z\in\Bbb Z^3}I(\alpha,\beta,\gamma)\phi_{\alpha,\beta,\gamma}(z).
$$
We now choose a Haar measure $\lambda$  on $H_3(\Bbb R)$ in such a way that $\lambda(I(\alpha,\beta,\gamma))=8\alpha\beta\gamma$ for some (and hence for all) positive reals $\alpha,\beta,\gamma$.

The proof of the following lemma is straightforward.

\proclaim{Lemma 3.1}
Let $\alpha_n,\beta_n,\gamma_n\to\infty$, $\alpha_n/
\gamma_n\to 0$ and
$\beta_n/\gamma_n\to 0$.
Then $(I(\alpha_n,\beta_n,\gamma_n))_{n=1}^\infty$
is a two-sided F{\o}lner sequence in $H_3(\Bbb R)$.
Moreover, if $(\alpha_n\beta_n)/\gamma_n\to 0$ then
$$
\frac{\lambda(I(\alpha_n,\beta_n,\gamma_n)\triangle I(\alpha_n,\beta_n,\gamma_n)^{-1})}{\lambda (I(\alpha_n,\beta_n,\gamma_n))}\to 0.
$$
\endproclaim

We now construct inductively the sequence $(C_{n+1},F_n)_{n=0}^\infty$.
Suppose that on the $n$-th step we already defined $F_0$, $(C_{j},F_j)_{j=1}^n$
and an auxiliary sequence $(\widetilde F_j)_{j=0}^{n-1}$ such that
$F_j=I(\alpha_j,\alpha_j,\gamma_j)$ for some $\alpha_j,\gamma_j>0$, $0\le j\le n$, and
$\widetilde F_j=I(\widetilde\alpha_j,\widetilde\alpha_j,\widetilde\gamma_j)$ for
for some $\widetilde\alpha_j,\widetilde\gamma_j>0$, $0\le j< n$.
Suppose also that  $F_j$ is equipped with a finite Borel partition
$\xi_j$, $0\le j<n$.
Our purpose is to construct $\xi_n,\widetilde F_n,  C_{n+1}$ and $F_{n+1}$.

Let $\xi_n$ be a finite Borel partition of  $F_n$  such that the following are satisfied:
\roster
\item"(i)"
the diameter (with respect to a natural metric on $h_3(\Bbb R)$) of each atom of $\xi_n$ is less than $\frac 1n$ and
\item"(iii)"
for each atom $A$ of  $\xi_{n-1}$ and each element $c\in C_{n}$,
the subset $Ac\subset F_{n}$ is $\xi_n$-measurable.
\endroster
Next, we introduce an auxiliary set
$$
S_n :=I((2n+1)\widetilde\alpha_{n-1},(2n+1)\widetilde\alpha_{n-1},
(2n+1)\widetilde\gamma_{n-1}).
$$
Then we select $\widetilde\alpha_n$ and $\widetilde\gamma_n$ to be the smallest positive reals such that
$$
F_nS_n\subset I(\widetilde\alpha_n,\widetilde\alpha_n,\widetilde\gamma_n).
$$
Now we set $\widetilde F_n:=I(\widetilde\alpha_n,\widetilde\alpha_n,\widetilde\gamma_n)$ and
 $\phi_n:=\phi_{\widetilde\alpha_n,\widetilde\alpha_n,\widetilde\gamma_n}$.
Next, for some integer $r_n>0$ (to be specified below), we let
$$
H_n:=\{(t_1,t_2,t_3)\in\Bbb Z^3\mid |t_1|<n^3,|t_2|<n^3, |t_3|<r_n\}.
$$
Suppose we are given a mapping $s_n:H_n\to S_n$.
Then we define another mapping $c_{n+1}:H_n\to H_3(\Bbb R)$ by setting
$c_{n+1}(h):=s_n(h)\phi_n(h)$.
We now set
$$
C_{n+1}:=c_{n+1}(H_n).
$$
Finally, let $F_{n+1}:=I(\alpha_{n+1},\alpha_{n+1},\gamma_{n+1})$,
where $\alpha_{n+1}$ and $\gamma_{n+1}$ are the smallest positive reals
such that $I(\alpha_{n+1},\alpha_{n+1},\gamma_{n+1})\supset  F_nC_{n+1}$.
It remains to specify $r_n$ and $s_n$.
For that we need an auxiliary lemma from \cite{dJ}.
In order to state it we first introduce some notation.
Given a  finite measure  $\nu$ on a finite
set $D$, we let $\|\nu\|_1:=\sum_{d\in D}|\nu(d)|$.
Given a finite set $Y$  and  a mapping $s:Y\to D$,
 let distr$_{y\in Y}s(y)$ denote the image of the equidistribution on $Y$ under $s$, i.e.
$$
(\text{distr}_{y\in Y}s(y))(d):=\frac {\# (s^{-1}(\{d\}))}{\# Y} \quad\text{for each}\quad d\in D.
$$

\proclaim{Lemma 3.2 \cite{dJ}}
Let $D$ be a finite set.
Then given $\epsilon >0$ and $\delta>0$,
 there  is  $R>0$ such that for each $r>R$,
there exists a map $s:\{-r,-r+1,\dots,r\}\to D$ such that
$$
\|\text{\rom{distr}}_{0\le t<N}(s(h+t),s(h'+t))-\lambda_{D}\times\lambda_D\|_1<\epsilon
$$
for each $N>\delta r$ and  $h\ne h'$ with $h+N<r$ and $h'+N<r$.
\endproclaim

For a finite subset $D$ in $S_n$, we denote by $\lambda_D$ the corresponding normalized {\it Dirac comb}, i.e. a measure  on $S_n$ given by
$\lambda_D(A):=\#(A\cap D)/\# D$ for each subset $A\subset S_n$.
Given two subsets $A,B\subset F_n$,
we let
$$
f_{A,B}(x,y):=\frac{\lambda(Ax\cap By)}{\lambda(F_n)}
$$
for all $x,y\in S_n$.
Now we choose a finite subset $D_n$ in $S_n$ such that
$$
\left|\int_{S_n\times S_n} f_{A,B}\,d\lambda_{D_n}d\lambda_{D_n}
-\frac 1{\lambda(S_n)^2}\int_{S_n\times S_n} f_{A,B}\,d\lambda d\lambda\right|\le \frac 1n
\tag3-1
$$
for all $\xi_n$-measurable subsets $A,B\subset F_n$ with $AS_n\cup BS_n\subset F_n$.
It follows from Lemma~3.2 that there exist $r_n>0$ and a mapping $s_n:H_n\to D_n$
such that $s_n(t_1,t_2,t_3)=s_n(0,0,t_3)$ for all $(t_1,t_2,t_3)\in H_n$ and
$$
\|\text{\rom{distr}}_{0\le t<N}(s_n(h+(0,0,t)),s_n(h'+(0,0,t)))- \lambda_{D_n}\times\lambda_{D_n}\|_1<\frac 1n\tag3-2
$$
for each $N>n^{-2}r_n$ and  $h,h'\in H_n$ with $h_3\ne h'_3$ and $h_3+N<r_n$ and $h'_3+N<r_n$.
Here $h_3$ and $h_3'$ denote the third coordinate of $h$ and $h'$ respectively.

Thus a sequence $(C_{n+1}, F_n)_{n=0}^\infty$ is well defined.
By Lemma~3.2  we may assume without loss of generality that $r_n\to\infty$ faster than exponentially.
Then the condition (I) from Section~2 is satisfied by Lemma~3.1.
It is straightforward to check that the conditions (II)--(IV)
and
 \thetag{2-2}  are also satisfied.
Hence the  $(C,F)$-action $T=(T_g)_{g\in H_3(\Bbb R)}$ of $H_3(\Bbb R)$
  associated with $(C_{n+1},F_n)_{n=0}^\infty$
is well defined on a standard non-atomic probability space $(X,\goth B,\mu)$.

We will need the following technical lemma.

\proclaim{Lemma 3.3}
Let $A,B$ and $S$ be subsets of finite Haar measure in $H_3(\Bbb R)$.
Then
$$
\int_{S\times S}\lambda(At_1\cap Bt_2)\,d\lambda(t_1)d\lambda(t_2)
=
\int_{A\times B}\lambda(aS\cap bS)\,d\lambda(a)d\lambda(b)
$$
\endproclaim
\demo{Proof}
Without loss of generality we may assume that
the subsets $A,B$ and $S$ are open and relatively compact in $H_3(\Bbb R)$.
Let
$$
\align
\Omega_1&:=\{(t_1,t_2,t_3)\in H_3(\Bbb R)^3\mid t_1,t_2\in S, t_3\in At_1\cap Bt_2\},\\
\Omega_2&:=\{(a,b,c)\in H_3(\Bbb R)^3\mid a\in A,b\in B, c\in aS\cap bS\}.
\endalign
$$
Then $\Omega_j$ is open and relatively compact subset in $H_3(\Bbb R)^3$, $j=1,2$.
Choosing an appropriate normalization of $\lambda$ we may assume without loss of generality that $\lambda^3(\Omega_1\cup \Omega_2)=1$.
For each $\epsilon>0$, we find  a lattice $\Gamma$ in $H_3(\Bbb R)$ such that
$$
\bigg|\frac{\#(\Gamma^3\cap \Omega_j)}{\#(\Gamma^3\cap(\Omega_1\cup \Omega_2))}-\lambda^3(\Omega_j)\bigg|<\epsilon, \quad j=1,2.
$$
Define a measure $\nu$ on $H_3(\Bbb R)$ by setting
$\nu(A):=\alpha\#(A\cap\Gamma)$ for all Borel subsets $A\subset H_3(\Bbb R)$, where the normalizing constant $\alpha>0$ is chosen in such a way
that $\nu^3(\Omega_1\cup \Omega_2)=1$.
Then
$$
|\nu^3(\Omega_j)-\lambda^3(\Omega_j)|<\epsilon, \quad j=1,2.\tag3-3
$$
We set $A_\Gamma:=A\cap\Gamma$, $B_\Gamma:=B\cap\Gamma$ and $S_\Gamma:=S\cap\Gamma$.
Then
$$
\aligned
\nu(At_1\cap Bt_2)&=\nu(A_\Gamma t_1\cap B_\Gamma t_2)\quad\text{and}\\
\nu(t_1A\cap t_2B)&=\nu(t_1A_\Gamma \cap t_2B_\Gamma)
\endaligned
\tag{3-4}
$$
whenever $t_1,t_2\in\Gamma$.
Applying Fubini theorem, \thetag{3-3} and \thetag{3-4} we obtain that
$$
\align
\int_{S\times S}\lambda(At_1\cap Bt_2)\,d\lambda(t_1)d\lambda(t_2)
&=\lambda^3(\Omega_1)\\
&=\nu^3(\Omega_1)\pm\epsilon\\
&=\int_{S_\Gamma\times S_\Gamma}\nu(At_1\cap Bt_2)\,d\nu(t_1)d\nu(t_2)\pm\epsilon\\
&=\int_{S_\Gamma\times S_\Gamma}\nu(A_\Gamma t_1\cap B_\Gamma t_2)\,d\nu(t_1)d\nu(t_2)\pm\epsilon\\
&=\sum_{a\in  A_\Gamma}\sum_{b\in B_\Gamma}\int_{S_\Gamma\times S_\Gamma}\nu(\{a t_1\}\cap \{b t_2\})\,d\nu(t_1)d\nu(t_2)\pm\epsilon\\
&=\int_{A_\Gamma\times B_\Gamma}\nu(a S_\Gamma\cap b S_\Gamma)\,d\nu(a)d\nu(b)\pm\epsilon\\
&=\int_{A\times B}\nu(a S\cap b S)\,d\nu(a)d\nu(b)\pm\epsilon\\
&=\nu^3(\Omega_2)\pm\epsilon\\
&=\lambda^3(\Omega_2)\pm 2\epsilon\\
&=\int_{A\times B}\lambda(aS\cap bS)\,d\lambda(a)d\lambda(b)\pm 2\epsilon.
\endalign
$$
\qed
\enddemo

We need some notation.
Given subsets $A_n,A_n'\subset F_n$, we write
$A_n\sim A_n'$ as $n\to\infty$ if $\lim_{n\to\infty}\lambda(A_n\triangle A_n')/\lambda(F_n)=0$.
This property implies that $\mu([A_n]_n\triangle[A_n']_n)\to 0$ as $n\to\infty$.
Let $e_3:=(0,0,1)\in\Bbb Z^3$.

\proclaim{Lemma 3.4}
$$
\sup_{A^*,B^*\subset F_{n-1}}|\mu(T_{\phi_n(e_3)}[A^*]_{n-1}\cap[B^*]_{n-1})-
\mu([A^*]_{n-1})\mu([B^*]_{n-1})|\to 0
$$
as $n\to\infty$.
More generally,
given a sequence of subsets $H_n^*\subset H_n$ such that
$$
\frac{\# (H_n^*\cap (H_n^*-e_3))}{\# H_n^*}\to 1\quad\text{and}\quad
\frac{\# H_n^*}{\# H_n}\to\delta
$$
for some $\delta>0$, we let $C_n^*:=\phi_n(H_n^*)$.
Then
$$
\sup_{A^*,B^*\subset F_{n-1}}|\mu(T_{\phi_n(e_3)}[A^*C_n^*]_{n}\cap[B^*]_{n-1})-
\mu([A^*C_{n}^*]_{n})\mu([B^*]_{n-1})|\to 0
$$
as $n\to\infty$.
\endproclaim

\demo{Proof}
Fix $n>0$ and two Borel subsets $A,B\subset F_n$.
Since $\phi_n(e_3)$ belongs to the center of $H_3(\Bbb R)$,
$$
\phi_n(e_3)c_{n+1}(h)=s_n(h)s_n(h+e_3)^{-1}c_{n+1}(h+e_3)\tag3-5
$$
for all $h\in\Bbb Z^3$ such that  $h, h+e_3\in H_n$.
We let
$$
F_n^{\circ}:=\{f\in F_n\mid fS_nS_n^{-1}\subset F_{n}\}.
$$
In other words, $F_n^\circ$ in an {\it interior} of $F_n$.
Of course, $F_n\sim F_n^\circ$ as $n\to\infty$.
We now put $A^\circ:=A\cap F_n^\circ$,
$B^\circ:=B\cap F_n^\circ$ and  $H_n(e_3):=H_n\cap(H_n-e_3)$.
Applying the standard approximation argument we may
assume without loss of generality that $A^\circ$ and $B^\circ$
are $\xi_n$-measurable.
It follows from~\thetag{2-4} and \thetag{2-5} that
$$
\aligned
\mu(T_{\phi_n(e_3)}[A]_n\cap[B]_n)&=\sum_{h\in H_n}\mu(T_{\phi_n(e_3)}[Ac_{n+1}(h)]_{n+1}\cap[B]_n)\\
&=\sum_{h\in H_n(e_3)}\mu([\phi_n(e_3)A^\circ c_{n+1}(h)]_{n+1}\cap[B]_n)+\overline o(1).
\endaligned
\tag3-6
$$
By $\overline o(1)$ here and below we mean a sequence that tends to $0$
uniformly in $A\subset F_n$ and $B\subset F_n$.
We now  deduce from~\thetag{3-5}, \thetag{3-6} and \thetag{2-3}, \thetag{2-4}, \thetag{2-6} that
$$
\align
\mu(T_{\phi_n(e_3)}[A]_n\cap[B]_n)
&=
\frac{1}{\# H_n}
\sum_{h\in H_n(e_3)}
\mu([(A^\circ s_n(h)s_n(h+e_3)^{-1}\cap B)]_{n})+\overline o(1).
\endalign
$$
Now \thetag{3-2}, \thetag{2-7}  and \thetag{3-1} yield
$$
\aligned
\mu(T_{\phi_n(e_3)}[A]_n\cap[B]_n)
&=
\frac{1}{(\# D_n)^2}
\sum_{d_1,d_2\in D_n}
\mu([(A^\circ d_1d_2^{-1}\cap B)]_{n})+\overline o(1)\\
&=
\frac{1}{(\# D_n)^2}
\sum_{d_1,d_2\in D_n}
\frac{\lambda(A^\circ d_1\cap B^\circ d_2)}{\lambda(F_n)}+\overline o(1)\\
&=\int_{S_n\times S_n}f_{A^\circ,B^\circ}\,d\lambda_{D_n}\lambda_{D_n}
+\overline o(1)\\
&=
\frac{1}{\lambda(S_n)^2}
\int_{S_n\times S_n}
f_{A^\circ,B^\circ}
\,d\lambda d\lambda+\overline o(1)\\
&=
\frac{1}{\lambda(S_n)^2}
\int_{S_n\times S_n}
f_{A,B}
\,d\lambda d\lambda+\overline o(1).
\endaligned
\tag3-7
$$
Suppose now that $A=A^*C_{n}$ and $B=B^*C_n$ for some subsets $A^*$ and $B^*$ in $F_{n-1}$.
We say that elements $c$ and $c'$ of $C_n$ are {\it partners} if $F_{n-1}cS_n\cap F_{n-1}c'S_n\ne\emptyset$.
We then write $c\bowtie c'$.
It follows that
$$
\int_{S_n\times S_n}
f_{A,B}
\,d\lambda d\lambda
=
\int_{S_n\times S_n}
\frac{\sum_{c\,\bowtie \,c'\in C_n}\lambda(A^*c t_1\cap B^* c' t_2)}{\lambda(F_{n})}\,d\lambda(t_1)d\lambda(t_2).
$$
Applying Lemma~3.3 we now obtain that
$$
\int_{S_n\times S_n}
f_{A,B}
\,d\lambda d\lambda
=\sum_{c\,\bowtie \,c'\in C_n}
\int_{A^*\times B^*}
\frac{\lambda(ac S_n\cap b c' S_n)}{\lambda(F_{n})}\,d\lambda(a)d\lambda(b).
\tag3-8
$$
Next, we note that
$$
\sup_{c,c'\in C_n}\sup_{a,b\in F_{n-1}}\bigg|\frac{\lambda(ac S_n\cap b c' S_n)-\lambda(c S_n\cap  c' S_n)}{\lambda(S_n)}\bigg|\to 0
\quad\text{as }n\to\infty.
$$
Hence it follows from \thetag{3-7} and \thetag{3-8} that
$$
\multline
\mu(T_{\phi_n(e_3)}[A^*]_{n-1}\cap[B^*]_{n-1})\\
=\frac{\sum_{c\,\bowtie\,c'\in C_n}}{\lambda(S_n)^2}
\int_{A^*\times B^*}
\frac{\lambda(c S_n\cap  c' S_n)+\lambda(S_n)\cdot\overline o(1)}{\lambda(F_{n})}\,d\lambda(a)d\lambda(b)+\overline o(1).
\endmultline
$$
A routine verification shows that
every
 $c\in C_n$ has no more than $2^{10}n^3$ partners.
Therefore
$$
\align
\mu(T_{\phi_n(e_3)}&[A^*]_{n-1}\cap[B^*]_{n-1})\\
&=\lambda(A^*)\lambda(B^*)\theta_n \pm
\frac{2^{10}n^3\#C_n}{\lambda(S_n)^2} \cdot \frac{\lambda(F_{n-1})^2\lambda(S_n)\cdot\overline o(1)}{\lambda(F_{n})}+\overline o(1)
\endalign
$$
for some $\theta_n>0$.
Since
$$
\lim_{n\to\infty}\frac{n^3\lambda(F_{n-1})}{\lambda(S_n)}>0
 \quad\text{and}\quad \lim_{n\to\infty}\frac{\lambda(F_{n-1})\#C_n}{\lambda(F_n)}=1,
$$
we obtain
$$
\mu(T_{\phi_n(e_3)}[A^*]_{n-1}\cap[B^*]_{n-1})
=\frac{\lambda(A^*)\lambda(B^*)}{\lambda(F_{n-1})^2}
\theta_n' +
\overline o(1)
$$
for some $\theta_n'>0$.
Substituting $A^*=B^*=F_{n-1}$ we obtain that $\theta_n\to 1$ as $n\to\infty$.
Thus
$$
\mu(T_{\phi_n(e_3)}[A^*]_{n-1}\cap[B^*]_{n-1})=\mu([A^*]_{n-1})
\mu([B^*]_{n-1})+
\overline o(1),
$$
and the first claim of the lemma is proved.
The second claim is proved in a similar way.
\qed
\enddemo

\proclaim{Corollary 3.5}
The sequence $(\phi_n(e_3))_{n=1}^{\infty}$
is mixing for $T$, i.e. $\mu(T_{\phi_n(e_3)}A\cap B)\to\mu(A)\mu(B)$ as $n\to\infty$ for all Borel subsets $A,B\subset X$.
Hence the flow $(T_{c(t)})_{t\in\Bbb R}$ is weakly mixing.
\endproclaim

We now refine this corollary.

\proclaim{Proposition 3.6}
The flow $(T_{c(t)})_{t\in\Bbb R}$ is mixing.
\endproclaim
\demo{Proof}
Take a sequence of reals  $t_n>0$ such that $g_n:=c(t_n)\in F_{n+1}\setminus F_n$, $n\in\Bbb N$.
It suffices to show that the sequence $(g_n)_{n=1}^\infty$ (or, at least a subsequence of it) is  mixing for $T$.
We write $g_n=f_n\phi_n(0,0,j_n)$ for some $f_n\in F_n\cap\{c(t)\mid t>0\}$ and $0<j_n<r_n$.
Let $h_n:=(0,0,j_n)\in\Bbb Z^3$, $H_n(h_n):=H_n\cap (H_n-h_n)$
and $F_n(f_n):=F_n\cap(f_n^{-1}F_n)$.
Passing to a subsequence, if necessary,  we may assume without loss of generality that
there exist $\delta_1\ge 0$ and $\delta_2\ge 0$ such that
$$
\frac{\# H_n(h_n)}{\# H_n}\to\delta_1\quad\text{and}\quad\frac{\lambda(F_n(f_n))}{\lambda(F_n) }\to\delta_2.
$$
Partition $C_{n+1}$ into three subsets $C^1_{n+1}$, $C^2_{n+1}$ and $C^3_{n+1}$ as follows
$$
\align
C^1_{n+1}&:=\{c\in C_{n+1}\mid g_nF_nc\subset F_{n+1}\phi_{n+1}(e_3)\},\\
C^2_{n+1}&:=\{c\in C_{n+1}\mid g_nF_nc\subset F_{n+1}\},\\
C^3_{n+1}&:=C_{n+1}\setminus (C^1_{n+1}\cup C^2_{n+1}).
\endalign
$$
We will show mixing separately on each of the subsets $[F_nC^1_{n+1}]_{n+1}$,
$[F_nC^2_{n+1}]_{n+1}$ and  $[F_nC^3_{n+1}]_{n+1}$ of $X$.

We first note that $\# C^3_{n+1}/\# C_{n+1}\to 0$ as $n\to\infty$.
Then \thetag{2-6} yields that
 $$
\lim_{n\to\infty}\sup_{A\subset F_n}\mu([AC^3_{n+1}]_{n+1})= 0.
\tag3-9
$$
Next, we note that $\phi_{n+1}(e_3)^{-1}g_nF_nC^1_{n+1}\subset F_{n+1}$ and hence
$$
T_{g_n}[AC^1_{n+1}]_{n+1}= T_{\phi_{n+1}(e_3)}[\phi_{n+1}(e_3)^{-1}g_nAC^1_{n+1}]_{n+1}.
$$
Hence, by the second claim of Lemma~3.4,
$$
\lim_{n\to\infty}\sup_{A,B\subset F_n}|\mu(T_{g_n}[AC^1_{n+1}]_{n+1}\cap[B]_n)
-\mu([AC^1_{n+1}]_{n+1})\mu([B]_n)|= 0.
\tag3-10
$$
It remains to consider the third case involving $C^2_{n+1}$.
It is a routine to verify that
$$
\frac{\# (C^2_{n+1}\triangle\{c_{n+1}(h)\mid h\in H_n(h_n)\}}
{\# C^2_{n+1}}\to 1,
$$
If $\delta_1=0$ then
$$
\lim_{n\to\infty}\sup_{A\subset F_n}\mu([AC^2_{n+1}]_{n+1})= 0.
\tag3-11
$$
Consider now the case where $\delta_1>0$.
Take two subsets
 $A,B\subset F_n$.
Partition $A$ into two subsets $A_1$ and $A_2$ such that
$f_nA_1\subset F_n$
and $f_nA_2\subset F_n\phi_n(e_3)$.
We note that $A_1=A\cap F_n(f_n)$.
Define $F_n^\circ$ in the same way as
 in the proof of Lemma~3.4
and set $A^\circ_1:=A_1\cap F_n^\circ$ and
$B^\circ:=B\cap F_n^\circ$.
Slightly modifying the reasoning in the proof of Lemma~3.4 we obtain
$$
\aligned
\mu(T_{g_n}[A_1C_{n+1}^2
&]_{n+1}
\cap[B]_n)=\sum_{h\in H_n(h_n)}\mu(T_{\phi_n(h_n)}[f_nA_1c_{n+1}(h)]_{n+1}\cap[B]_n)
+\overline o(1)\\
&=\sum_{h\in H_n(h_n)}\mu([\phi_n(t_n)f_nA_1^\circ c_{n+1}(h)]_{n+1}\cap[B]_n)+\overline o(1)\\
&=
\frac{\delta_1}{\# H_n(t_n)}
\sum_{h\in H_n(t_n)}
\mu([(f_nA_1^\circ s_n(h)s_n(h+t_n)^{-1}\cap B^\circ)]_{n})+\overline o(1)\\
&=
\frac{\delta_1}{\lambda(S_n)^2}
\int_{S_n\times S_n}
f_{A_1,B}
\,d\lambda d\lambda+\overline o(1).
\endaligned
$$
If $\delta_2=0$ then
$$
\lim_{n\to\infty}\sup_{A\subset F_n}\mu([A_1C^2_{n+1}]_{n+1})= 0.
\tag3-12
$$
Consider now the case where $\delta_2>0$.
As in the proof of Lemma~3.4 we now take $A=A^*C_n$ and $B=B^*C_n$ for some Borel subsets $A^*,B^*\subset F_{n-1}$.
Let $C_n':=C_n\cap F_n(f_n)$.
It follows that
$$
\frac{\# C_n'}{\# C_n}\to\delta_2\quad\text{and}\quad
\sup_{A^*\subset F_{n-1}}\mu([A_1]_n\triangle[A^*C_n']_n)\to 0
$$
as $n\to\infty$.
Hence $\sup_{A^*\subset F_{n-1}}|\mu([A_1]_n)-\delta_2\mu([A^*]_{n-1})|\to 0$.
Arguing as in the proof of Lemma~3.4 we obtain that
$$
\sup_{A^*,B^*\subset F_{n-1}}|\mu(T_{g_n}[A^* C_n'C_{n+1}^2]_{n+1}\cap[B^*]_{n-1})-
\delta_2\mu([A^*]_{n-1})
\mu([B^*]_{n-1})|\to 0.
$$
Therefore
$$
\lim_{n\to\infty}\sup_{A^*,B^*\subset F_{n-1}}|\mu(T_{g_n}[A_1]_{n}\cap[B^*]_{n-1})-
\mu([A_1]_{n})
\mu([B^*]_{n-1})|= 0.
\tag3-13
$$
Since
$T_{g_n}[A_2]_n=T_{\phi_n(h_n+e_3)}[\phi_n(e_3)^{-1}f_nA_2]_n$
with $\phi_n(e_3)^{-1}f_nA_2\subset F_n$,
a similar reasoning yields
$$
\lim_{n\to\infty}\sup_{A^*,B^*\subset F_{n-1}}|\mu(T_{g_n}[A_2]_{n}\cap[B^*]_{n-1})-
\mu([A_2]_{n})
\mu([B^*]_{n-1})|= 0.
\tag3-14
$$
Since
$$
[A^*]_{n-1}=[A^*C_nC_{n+1}^1]_{n+1}\sqcup [A^*C_nC_{n+1}^3]_{n+1}
\sqcup [A_1C_{n+1}^2]_{n+1}\sqcup [A_2C_{n+1}^2]_{n+1},
$$
it follows from \thetag{3-9}--\thetag{3-14} that
$$
\lim_{n\to\infty}\sup_{A^*,B^*\subset F_{n-1}}|\mu(T_{g_n}[A^*]_{n-1}\cap[B^*]_{n-1})-
\mu([A^*]_{n-1})
\mu([B^*]_{n-1})|= 0,
$$
 i.e.  $(g_n)_{n=1}^\infty$ is a mixing sequence for $T$, as desired.
\qed
\enddemo

We now state and prove the main result of this section.
By an advice of V.~Ryzhikov we deduce it from Proposition~3.6 by adapting the argument used in \cite{Ry1, Theorem 6} and \cite{Ry2, Theorem~4.4} to show mixing of $(T_{a(t)})_{t\in\Bbb R}$.

\proclaim{Theorem 3.7} Let $T$ be an action of $H_3(\Bbb R)$
such that the flow $(T_{c(t)})_{t\in\Bbb R}$ is ergodic.
Then $T$ is mixing.
\endproclaim
\demo{Proof}
If $T$ were not mixing then there exist $\epsilon>0$, a sequence $g_n\to\infty$ in $H_3(\Bbb R)$ and
two subsets $A_0,B_0
\subset X$ such that
 $$
|\mu(A_0\cap T_{g_n}B_0)-\mu(A_0)\mu(B_0)|>\epsilon\qquad\text{for all $n$.}\tag3-15
$$
We write $g_n$ as $g_n=a(t_{n,1})b(t_{n,2})c(t_{n,3})$ with $t_{n,1},t_{n,2},t_{n,3}\in\Bbb R$.
Since the flow $(T_{c(t)})_{t\in\Bbb R}$ is mixing by Proposition~3.6, we may assume without loss of generality (passing to a subsequence, if necessary) that either $t_{n,1}\to\infty$ or $t_{n,2}\to\infty$.
We consider the former case (the latter case is analogous).
Select $\gamma>0$ such that
$$
\sup_{n>0}|\mu(A_0\cap T_{g_n}B_0)-\mu(T_{b(\beta)}A_0\cap T_{g_n}T_{b(\beta)}B_0)|<\epsilon/2
\qquad\text{for all $\beta<\gamma$.}
$$
Consider a probability measure $\kappa_{\gamma,n}$ on $X\times X$ by setting
$$
\kappa_{\gamma,n}(A\times B):=\frac 1\gamma\int_0^\gamma\mu(A\cap T_{b(\beta)^{-1}g_nb(\beta)}B)\,d\beta \qquad\text{for all $A,B\subset X$.}
$$
It is easy to see that $\kappa_{\gamma,n}$ is a 2-fold self-joining of $(T_{c(t)})_{t\in\Bbb R}$
and
$$
|\kappa_{\gamma,n}(A_0\times B_0)-\mu(A_0\cap T_{g_n}B_0)|<\epsilon/2\qquad
\text{for all $n$.}
\tag{3-16}
$$
Passing to a further subsequence we may assume without loss of generality that the sequence $(\kappa_{\gamma,n})_{n\in\Bbb N}$ converges weakly to a 2-fold self-joining $\kappa_\gamma$ of $(T_{c(t)})_{t\in\Bbb R}$.
Since
$$
{b(\beta)^{-1}g_nb(\beta)}={g_nc(\beta t_{n,1})} \qquad\text{for all $\beta\in\Bbb R$},
$$
 it follows that
$$
\align
\kappa_{\gamma,n}(A\times T_{c(t)}B)
&=\frac 1\gamma\int_{0}^\gamma\mu(A\cap T_{g_nc(\beta t_{n,1}+t)}B)\,d\beta\\
&=\frac 1\gamma\int_{t/t_{n,1}}^{\gamma+t/t_{n,1}}\mu(A\cap T_{g_nc(\beta t_{n,1})}B)\,d\beta.
\endalign
$$
Hence $\kappa_{\gamma,n}(A\times T_{c(t)}B)
\to \kappa_{\gamma}(A\times B)$ as $n\to\infty$.
This means that $\kappa_\gamma\circ(\text{Id}\times T_{c(t)})=\kappa_\gamma$ for all $t\in\Bbb R$.
This yields that $\kappa_\gamma=\mu\times\mu$.
We obtain a contradiction with \thetag{3-15} plus \thetag{3-16}.
\qed
\enddemo

We recall that the group Aut$(X,\mu)$ of all invertible $\mu$-preserving transformations of $X$ is Polish with respect to the weak topology.
 The {\it weak topology} is the weakest topology in which all the maps Aut$(X,\mu)\ni S\mapsto\mu(SA\cap B)\in\Bbb R$, $A,B\subset X$, are continuous.

As a byproduct of the proof of Theorem~3.7, we obtain the following results.

\proclaim{Corollary 3.8}
Let $T$ be an action of $H_3(\Bbb R)$
such that the flow $(T_{c(t)})_{t\in\Bbb R}$ is ergodic.
Then the following are satisfied.
\roster
\item"(i)"
Every sequence $(g_n)_{n=1}^\infty$ in $H_3(\Bbb R)$
such that $g_n=a(t_{n,1})b(t_{n,2})c(t_{n,3})$ with $(t_{1,n})_{n=1}^\infty$ or $(t_{n,2})_{n=1}^\infty$ unbounded is mixing.
\item"(ii)"
The flows  $(T_{a(t)})_{t\in\Bbb R}$ and
$(T_{b(t)})_{t\in\Bbb R}$ are both mixing \cite{Ry1, Theorem 6}
 and \cite{Ry2, Theorem~4.4}).
\item"(iii)"
The weak closure of $\{T_g\mid g\in H_3(\Bbb R)\}$ in \rom{Aut}$(X,\mu)$ is the union of
$\{T_g\mid g\in H_3(\Bbb R)\}$ and the weak closure of
$\{T_{c(t)}\mid t\in \Bbb R\}$.
\item"(iv)"
If $T$ is rigid then $\{T_{c(t)}\mid t\in \Bbb R\}$ is rigid.
\endroster
\endproclaim

\remark{Remark \rom{3.9}} In a similar way, we can show the following.
Let $T$ be an action of $H_3(\Bbb R)$.
\roster
\item"(i)"
If
 the flow $(T_{c(t)})_{t\in\Bbb R}$ is
 mixing of order $k$
 then $T$ is mixing of order $k$.
\item"(ii)"
If $(T_{c(t)})_{t\in \Bbb R}$ is weakly mixing then
the flows
$(T_{a(t)})_{t\in\Bbb R}$ and
$(T_{b(t)})_{t\in\Bbb R}$ are both
 mixing of order $k$
\cite{Ry1, Theorem 6}
 and \cite{Ry2, Theorem~4.4}).
\endroster
Then a natural question arises: how can we  construct $T$ such that
$(T_{c(t)})_{t\in\Bbb R}$ is
 mixing of order $k$
for $k>1$?
For that we  need to modify slightly the construction of $T$ from Proposition~3.6.
Indeed, the  only modification concerns the right choice of the `spacer mappings' $s_n$.
Namely, we need to replace \thetag{3-2} with the following subtler condition:
$$
\|\text{\rom{distr}}_{0\le t<N}(s_n(h^{(1)}+(0,0,t)),\dots,s_n(h^{(k)}+(0,0,t)))- \lambda_{D_n}^k\|_1<\frac 1n
$$
for each $N>n^{-2}r_n$ and  $h^{(1)},\dots, h^{(k)}\in H_n$ with $h^{(1)}_3\ne h^{(2)}_3\ne \cdots\ne h_3^{(k)}$ and $h_3^{(1)}+N<r_n$,\dots, $h^{(k)}_3+N<r_n$.
We leave details to the readers.
\endremark

We conclude this section with a simple but useful observation.

\proclaim{Proposition 3.10}
Let $T$ be an action of $H_3(\Bbb R)$.
If $(T_{c(t)})_{t\in\Bbb R}$
is ergodic then it is weakly mixing.
\endproclaim
\demo{Proof}
Indeed, let $f\circ T_{c(t)}=e^{ist}f$ for some $0\ne s
\in\Bbb R$ and $0\ne f\in L^2(X,\mu)$.
Then for each $g\in H_3(\Bbb R)$, $(f\circ T_g)\circ T_{c(t)}=e^{ist}f\circ T_g$.
Since the center is ergodic, it follows that $f\circ T_g=\xi(g)f$ for some $\xi(g)\in\Bbb C$.
It is straightforward that $\xi$ is a continuous character of $H_3(\Bbb R)$.
Since every continuous character of $H_3(\Bbb R)$ is trivial on the center, it follows that $f$ is invariant under $(T_{c(t)})_{t\in\Bbb R}$.
Hence it is constant.
\qed
\enddemo

\head 4. The `joining' structure of the action
\endhead

By a polyhedron in $\Bbb R^3$ we mean a union of finitely many mutually disjoint convex polyhedrons.
We say that a  polyhedron  is {\it rational} if  the  coordinates of its vertices  are all rational.
If no one face of a polyhedron is parallel  to the line $\{0\}\times\{0\}\times\Bbb R\subset\Bbb R^3$
then we call  this polyhedron  {\it normal}.
The intersection of two normal polyhedrons  and the union of two disjoint normal polyhedrons is a normal polyhedron.
It is a routine to verify that if $A$ is a normal polyhedron and $g\in H_3(\Bbb R)$ then
$gA$ and $Ag$ are normal polyhedrons.
 Given a  normal polyhedron $A$ and $\epsilon>0$, there is a rational polyhedron $B\subset A$ such that
$$
\sup_{t_1,t_2}\lambda_{\Bbb R}(\{t\mid (t_1,t_2,t)\in A\setminus B\})<\epsilon.
$$

From now on all the polyhedrons are assumed to be open subsets of $\Bbb R^3$.

\proclaim{Theorem 4.1} Let $T$ be the action of $H_3(\Bbb R)$ constructed in Section~3.
Then the following are satisfied.
\roster
\item"\rom{(i)}"
The flow $(T_{c(t)})_{t\in\Bbb R}$ is simple
and $C((T_{c(t)})_{t\in\Bbb R})=\{T_g\mid g\in H_3(\Bbb R)\}$.
\item"\rom{(ii)}"
The transformation $T_{c(1)}$ is simple
and $C(T_{c(1)})=\{T_g\mid g\in H_3(\Bbb R)\}$.
\endroster
\endproclaim

\demo{Proof}
(i) Take an ergodic 2-fold self-joining $\nu$ for $(T_{c(t)})_{t\in\Bbb R}$.
Identifying $H_3(\Bbb R)$ with $\Bbb R^3$ via the mapping $c(t_3)b(t_2)a(t_1)\mapsto (t_1,t_2,t_3)$ we can say about polyhedrons in $H_3(\Bbb R)$.
Recall that $F_n=I(\alpha_n,\alpha_n,\gamma_n)$.
We call a point $(x,x')\in X\times X$  {\it generic} for $(\nu,(T_{c(t)}\times T_{c(t)})_{t\in\Bbb R})$ if
$$
\frac {n^2}{\gamma_n}\int_0^{\gamma_n/n^2}
1_{[A]_m\times [A']_m}(T_{c(t)}x,T_{c(t)}x')\,d\lambda_\Bbb R(t)\to\nu([A]_m\times [A']_m)\tag4-1
$$
for each pair of rational polyhedrons $A,A'\subset F_m$, $m\in\Bbb N$.
It follows from  the pointwise ergodic theorem for ergodic flows that $\nu$-almost every point of $X\times X$ is generic for $(\nu,(T_{c(t)}\times T_{c(t)})_{t\in\Bbb R})$.
Fix such a point  $(x,x')$.
We now claim that \thetag{4-1} holds for each pair of normal convex (not necessarily rational) polyhedrons $A,A'\subset F_m$.
For that approximate $A$ and $A'$ with nested sequences of rational convex polyhedrons $A_1\subset A_2\subset\cdots\subset A$ and
$A_1'\subset A_2'\subset\cdots\subset A'$  such that
$$
\sup_{t_1,t_2}\lambda_{\Bbb R}(\{t\mid (t_1,t_2,t)\in (A\setminus A_j)\cup (A'\setminus A_j')\})<\frac 1j
\quad\text{for all $j$}
$$
and pass to the limit in \thetag{4-1}.

Since $\nu$ is a 2-fold self-joining of $T$, it follows that
$x$ and $x'$ are generic points for  $(\mu,T_{c(t)})_{t\in\Bbb R})$, i.e.
$$
\frac {n^2}{\gamma_n}\int_{0}^{\gamma_n/n^2}
1_{[A]_m}(T_{c(t)}x)\,d\mu_\Bbb R(t)\to\mu([A]_m)
$$
for each normal polyhedron $A$ in $F_m$, $m\in\Bbb N$.
If $n$ is sufficiently large
then we can write $x$ and $x'$ as infinite sequences $x=(f_n,c_{n+1}, c_{n+2},\dots)$ and $x'=(f_n',c_{n+1}', c_{n+2}',\dots)$
with $f_n,f_n'\in F_n$ and $c_j,c_j'\in C_j$ for all $j>n$.
Represent $f_n,f_n',c_j$ and $c_j'$ as products
$$
\align
f_n=c(\tau_{3,n})b(\tau_{2,n})a(\tau_{1,n}),\quad
&f_n'=c(\tau_{3,n}')b(\tau_{2,n}')a(\tau_{1,n}'),\\
c_j=c(t_{3,j})b(t_{2,j})a(t_{1,j}),\quad
&c_j'=c(t_{3,j}')b(t_{2,j}')a(t_{1,j}').
\endalign
$$
Since $F_n=I(\alpha_n,\alpha_n,\gamma_n)$, we obtain that
$|\tau_{3,n}|<\gamma_n$ and $|\tau_{3,n}'|<\gamma_n$.
Moreover, by a standard application of Borel-Cantelli lemma, we may assume without loss of generality that
$\frac{\max(t_{3,j},t'_{3,j})}{\gamma_j}<1-\frac{4}{j^2}$ for all $j>n$ if $n$
is sufficiently large.
Increasing $n$ by $1$ if necessary, we will assume also that
$\frac{\max(\tau_{3,n},\tau'_{3,n})}{\gamma_n}<1-\frac{1}{n^2}$.
We set $t_n:=f_n'f_n^{-1}$.

Let $B$ and $B'$ be two normal polyhedrons in $F_m$.
It follows from \thetag{4-1} that for each $\epsilon>0$ there is $n>m$
such that
$$
\frac {n^2}{\gamma_n}\int_{0}^{\gamma_n/n^2}
1_{[B]_m\times [B']_m}(T_{c(t)}x,T_{c(t)}x')
\,d\lambda_\Bbb R(t)=\nu([B]_m\times [B']_m)
\pm\epsilon.\tag4-2
$$
Let $\widetilde B:=BC_{m+1}\cdots C_n$
and $\widetilde B':=B'C_{m+1}\cdots C_n$.
Then $\widetilde B$ and $\widetilde B'$ are normal polyhedrons
in $F_n$ and $[B]_m=[\widetilde B]_n$ and
$[B']_m=[\widetilde B']_n$.
If $g\in H_3(\Bbb R)$ and $gf_n'\in F_n$
then $T_gx'\in [\widetilde B']_n$ if and only if $T_gT_{t_{n}}x\in[\widetilde B']_n$.
Hence
$$
1_{[\widetilde B]_n\times[\widetilde B']_n}(T_gx,T_gx')=1_{T_g^{-1}[\widetilde B]_n\cap T_{t_n}^{-1}T_g^{-1}[\widetilde B']_n}(x).
\tag4-3
$$
Since $n$ is large, $c(\gamma_n/n^2)f_n'\in F_n$ and therefore \thetag{4-2} and \thetag{4-3} yield
$$
\frac {n^2}{\gamma_n}\int_{0}^{\gamma_n/n^2}
1_{T_{c(t)}^{-1}([\widetilde B]_n\cap T_{t_n}^{-1}[\widetilde B']_n)}(x)
\,d\lambda_\Bbb R(t)=\nu([B]_m\times [B']_m)
\pm\epsilon.\tag4-4
$$
We note that for each $g\in H_3(\Bbb R)$ there is $M>n$ such that
$gF_nC_{n+1}\cdots C_M\subset F_M$.
Therefore choosing a sufficiently large $M$ we can write for each $0<t<\gamma_n/n^2$ that
$$
T_{c(t)}^{-1}([\widetilde B]_n\cap T_{t_n}^{-1}[\widetilde B']_n)
=[D]_M,
$$
where $D=c(t)^{-1}(\widetilde BC_{n+1}\cdots C_M\cap t_n^{-1}
\widetilde B'C_{n+1}\cdots C_M)$.
Since $D$ is a normal polyhedron in $F_M$ and
 $x$ is generic for $(\mu,T_{c(t)})_{t\in\Bbb R})$, we
deduce from \thetag{4-4} that
$$
\frac {n^2}{\gamma_n}\int_{0}^{\gamma_n/n^2}
\mu(T_{c(t)}^{-1}([\widetilde B]_n\cap T_{t_n}^{-1}[\widetilde B']_n))
\,d\lambda_\Bbb R(t)=\nu([B]_m\times [B']_m)
\pm\epsilon.
$$
Thus
$$
\mu([B]_m\cap T_{t_n}^{-1}[B']_m)= \nu([B]_m\times [B']_m)\pm\epsilon.
$$
We note that the sequence $(t_n)_{n=1}^\infty$  either converges to some $g\in H_3(\Bbb R)$ (in the case when $c_j=c_j'$ for all sufficiently large $j$) or tends to infinity (in the case when $c_j\ne c_j'$ for infinitely many $j$).
In the first case $\nu=\mu_{T_g}$ and in the second case $\nu=\mu\times\mu$.
In the second case we use the fact that $T$ is mixing.
Thus the flow
$(T_{c(t)})_{t\in\Bbb R}$ is 2-fold simple
and its centralizer
is
$\{T_g\mid g\in H_3(\Bbb R)\}$.
By \cite{Ry4}, each 2-fold simple flow is simple.

(ii) follows from (i) and \cite{dJR, Theorem~6.1}.
\qed
\enddemo

\remark{Remark 4.2}
In a similar way one can show a more general fact: each mixing $(C,F)$-action of $H_3(\Bbb R)$
with $F_n=I(\alpha_n,\beta_n,\gamma_n)$ and $\gamma_n\gg \alpha_n\beta_n$
is 2-fold simple.
\endremark

\proclaim{Theorem 4.3} Given an action $T$ of $H_3(\Bbb R)$,
let
the flow $(T_{c(t)})_{t\in\Bbb R}$ be simple
and $C((T_{c(t)})_{t\in\Bbb R})=\{T_g\mid g\in H_3(\Bbb R)\}$.
Then
\roster
\item"\rom{(i)}"
$T$  has MSJ and
\item"\rom{(ii)}"
the actions $(T_g)_{g\in H_{2,a}}$ and $(T_g)_{g\in H_{2,b}}$ have MSJ.
\endroster
\endproclaim

\demo{Proof}
(i) Let $\nu\in J_2^e(T)$.
It follows from Theorem~4.1(i) that there is a probability measure $\kappa$ on $H_3(\Bbb R)$ and a non-negative real $\theta\le 1$ such that
$$
\nu=\theta \mu\times\mu+(1-\theta)\int_{H_3(\Bbb R)}\mu_{T_g}\,d\kappa(g).
$$
Since $\nu$ is $(T\times T)$-ergodic and $\mu\times\mu\in J_2(T)$,
either $\theta=0$ or $\theta=1$.
We consider the former case.
 We note that $\nu$ is $(T\times T)$-invariant if and only if  $\kappa$ is invariant under all inner automorphisms of $H_3(\Bbb R)$.
Since $\kappa$ is finite, it is supported on the center of $H_3(\Bbb R)$.
The ergodicity of $\nu$ implies that $\kappa$ is a singleton.
Therefore $\nu=\mu_{T_{c(t)}}$ for some $t\in\Bbb R$.
Thus, $T$ has MSJ$_2$.
Now Theorem~4.1(i) implies that $T$ has MSJ.

(ii) is shown in a similar  way.
\qed
\enddemo

\head 5. Mixing Poisson and mixing Gaussian  actions of Heisenberg group
\endhead

We first construct a mixing infinite measure preserving rank-one action of $H_3(\Bbb R)$.
Let $(F_n,C_{n+1})_{n=0}^\infty$ be  sequence of subsets in $H_3(\Bbb R)$ satisfying (I)--(IV) from Section~2.
Suppose, in addition, that \thetag{2-2} is not satisfied.
Let $T$ be the $(C,F)$-action associated with
$(F_n,C_{n+1})_{n=0}^\infty$.
Let $(X,\goth B,\mu)$ denote the space of this action.
Then
$T$ is  of rank one along $(F_n)_{n=1}^\infty$ and $\mu(X)=\infty$.

\proclaim{Theorem 5.1}
Suppose that the following conditions are satisfied:
\roster
\item"(i)"
$F_nF_n^{-1}F_nC_n\subset F_{n+1}$,
\item"(ii)"
the sets
$F_nc_1c_2^{-1}F_n^{-1}$, $c_1\ne c_2\in C_{n+1}$, and $F_nF_n^{-1}$ are all pairwise disjoint and
\item"(iii)"
$\# C_n\to\infty$ as $n\to\infty$.
\endroster
 Then $T$ is mixing.
\endproclaim
\demo{Proof}
Let $(g_n)_{n=1}^\infty$ be a sequence in $H_3(\Bbb R)$ such that
$g_n\to\infty$.
We verify that this sequence (or, a subsequence of it) is mixing, i.e. $\mu(T_{g_n}D_1\cap D_2)\to 0$ for all subsets $D_1,D_2\subset X$ of finite measure.
Without loss of generality we may assume that $g_n\in F_nF_n^{-1}\setminus F_{n-1}F_{n-1}^{-1}$.
Take $A,B\subset F_n$.
It follows from~(i) and \thetag{2-3}--\thetag{2-6} that
$$
\align
\mu(T_{g_n}[A]_n\cap[B]_n)&=\mu([g_nAC_{n+1}\cap BC_{n+1}]_{n+1})\\
&=
\sum_{c_1,c_2\in C_{n+1}}\mu([g_nAc_{1}\cap Bc_{2}]_{n+1})\\
&=
\mu([g_nA\cap B]_n)+
\sum_{c_1\ne c_2\in C_{n+1}}\mu([g_nAc_{1}\cap Bc_{2}]_{n+1})
\endalign
$$
If $\mu([g_nAc_{1}\cap Bc_{2}]_{n+1})>0$ for some $c_1\ne c_2\in C_{n+1}$ then $g_n\in F_nc_2c_1^{-1}F_n^{-1}$.
Since $g_n\in F_nF_n^{-1}$, we obtain a contradiction with (ii).
Hence
$$
\mu(T_{g_n}[A]_n\cap[B]_n)=\mu([g_nA\cap B]_n).
$$
Suppose now that $A=A^*C_{n-1}$, $B=B^*C_{n-1}$ for some
subsets $A^*,B^*\subset F_{n-1}$.
Then
$$
\mu([g_nA\cap B]_n)=\sum_{c_1,c_2\in C_{n}}\mu([g_nA^*c_1\cap B^*c_2]_n).
$$
 It follows from (ii) that there is no more than one pair $(c_1,c_2)\in C_n\times C_n$ such that $c_1\ne c_2$ and $\mu([g_nA^*c_1\cap B^*c_2]_n)>0$.
Hence
$$
\mu([g_nA\cap B]_n)=\mu([g_nA^*\cap B^*]_{n-1})+\mu([g_nA^*c_1\cap B^*c_2]_n).
$$
If $\mu([g_nA^*\cap B^*]_{n-1})>0$ then
 $g_n\in F_{n-1}F_{n-1}^{-1}$, a contradiction.
On the other hand,
$$
\mu([g_nA^*c_1\cap B^*c_2]_n)\le
\mu([g_nA^*c_1]_{n})=\mu([A^*]_{n-1})/\# C_{n}
$$
by \thetag{2-6}.
Therefore
$$
\mu(T_{g_n}[A]_n\cap[B]_n)\le\mu([A^*]_{n-1})/\# C_n=
\mu([A]_{n})/\# C_n.
$$
It remains to use the standard approximation of $D_1$ and $D_2$ with cylinders $[A]_n$ and $[B]_n$ respectively and  apply (iii).
\qed
\enddemo

We now recall that given a $\sigma$-finite infinite non-atomic standard measure space $(X,\goth B,\mu)$, there is a canonical way to associate s standard  probability space  $(\widetilde X,\widetilde B,\widetilde\mu)$ which is called the {\it Poisson suspension} of $(X,\goth B,\mu)$ (see \cite{CFS}, \cite{Ro}).
Moreover, there exists a continuous homomorphism from the group Aut$(X,\mu)$ of $\mu$-preserving transformations of $X$ to the group Aut$(\widetilde X,\widetilde\mu)$ of $\widetilde\mu$-preserving transformations of $\widetilde X$.
The image of a transformation $S\in\text{Aut}(X,\mu)$ under this homomorphism  is denoted by $\widetilde S$.
It is called the {\it Poisson suspension} of $S$.
Given a $\mu$-preserving action $T$ of $H_3(\Bbb R)$, we consider
a $\widetilde\mu$-preserving action $\widetilde T=(\widetilde T_g)_{g\in H_3(\Bbb R)}$ of $H_3(\Bbb R)$ and call it the {\it Poisson suspension} of $T$.
Consider the Koopman representation $U_{T}$ of $H_3(\Bbb R)$ in $L^2(X,\mu)$.
Then we can associate to $U_T$ a (probability preserving) Gaussian
action $T^*=(T^*_g)_{g\in H_3(\Bbb R)}$ of $H_3(\Bbb R)$ (see e.g. \cite{Gl}).
The Koopman representations $U_{\widetilde T}$ and $U_{T^*}$ of $H_3(\Bbb R)$ are unitarily equivalent.

\proclaim{Corollary 5.2} Let $T$ be an action of $H_3(\Bbb R)$ constructed in Theorem~5.1.
Then the Poisson suspension   $\widetilde T$ of $T$ is mixing of all orders.
The corresponding  Gaussian action $T^*$ of $H_3(\Bbb R)$ is also mixing of all orders.
\endproclaim
\demo{Proof}
The first assertion follows from  \cite{Ro, Theorem~4.8}.
Since $U_{\widetilde T}$ and $U_{T^*}$ are unitarily equivalent
and the property of mixing is spectral (i.e. it is a property of the associated Koopman representation), we obtain that
 $T^*$ is also mixing.
The multiple mixing of $T^*$ follows now from \cite{Le}.
\qed
\enddemo

\head 6. Asymmetry in  actions of Heisenberg group
\endhead

In this section we construct  a $(C,F)$-action $T$ of $H_3(\Bbb R)$ with ergodic flow $(T_{c(t)})_{t\in\Bbb R}$
such that the transformation $T_{c(1)}$ is not conjugate to
$T_{c(1)}^{-1}$.

Suppose that $n$ is divisible by 3 and on the $n$-th step of the inductive construction of the sequence $(C_{n+1}, F_n)_{n=0}^\infty$
we have defined $F_n=I(\alpha_n,\alpha_n,\gamma_n)$ for some parameters
$\alpha_n$ and $\beta_n$.
We are going to define $C_{n+1}$ and $F_{n+1}$.
Let $\phi_n:=\phi_{\alpha_n,\alpha_n,\gamma_n}$.
We set
$$
s_n(j):=
\cases
c(0)& \text{if }j\equiv 0 \pmod 5,\\
c(1)& \text{if }j\equiv 1\pmod 5 \text{ or }j\equiv 2\pmod 5,\\
c(2)& \text{if }j\equiv 3\pmod 5 \text{ or }j\equiv 4\pmod 5,
\endcases
$$
 $c_{n+1}(j):=s_n(j)\phi_n(0,0,j)$ and $C_{n+1}:=\{c_{n+1}(j)\mid |j|\le n\}$.
Next, let $F_{n+1}:=I(\alpha_{n},\alpha_n,\gamma_{n+1})$, where $\gamma_{n+1}$ is the minimal positive real such that $F_{n}C_{n+1}\subset I(\alpha_{n},\alpha_n,\gamma_{n+1})$.

If $n$ is not divisible by 3,
we do the $n$-th construction step exactly  as in Section~3.

In such a way we define completely the sequence $(C_{n+1}, F_n)_{n=0}^\infty$.
The conditions (I)--(IV) and~\thetag{2-2} from Section~2 are all satisfied.
Denote by $T$ the associated   $(C,F)$-action of $H_3(\Bbb R)$.
The probability space of this action is denoted by $(X,\mu)$.

\proclaim{Theorem 6.1}
The flow $(T_{c(t)})_{t\in\Bbb R}$ is weakly mixing and rigid.
It is not conjugate to the flow $(T^{-1}_{c(t)})_{t\in\Bbb R}$.
If  $\gamma_n\in\Bbb N$ for all $n\in\Bbb N$ then
the transformation $T_{c(1)}$ is not conjugate to
$T_{c(1)}^{-1}$.
\endproclaim
\demo{Proof}
The fact that $(T_{c(t)})_{t\in\Bbb R}$ is weakly mixing follows immediately from the definition of $C_{n+1},F_{n+1}$ when  $3\nmid n$ and the proof of Lemma~3.4 (see also Corollary~3.5).
It follows from the definition of $C_{n+1}$ when $3\mid n$ that this flow
is rigid.
Indeed, it is easy to verify that $T_{\phi_n(0,0,5)}\to\text{Id}$ as $n\to\infty$ and $3 \mid n$.

We now prove the second claim.
Fix $n$ which is divisible by $3$.
Take subsets $A,B,C,D\subset X$.
Partition $C_{n+1}$ into $C_{n+1}^j$, $j=0,\dots,4$, as follows
$$
C^j_{n+1}:=\bigsqcup_{t\equiv j\pmod 5}F_nc_{n+1}(t)
$$
and set
$F^j_{n+1}:=F_nC_{n+1}^j$.
Let $l_n:=c(1)\phi_n(0,0,1)$.
We now claim that
$$
\align
\mu(F^0_{n+1}\cap A\cap T_{l_n}B\cap T_{l_n^2}C\cap T_{l_n^3}D) &\to\frac 15\mu(A\cap T_{c(1)}^{-1}B\cap T_{c(1)}^{-2}C\cap T_{c(1)}^{-2}D),
\\
\mu(F^1_{n+1}\cap A\cap T_{l_n}B\cap T_{l_n^2}C\cap T_{l_n^3}D) &\to\frac 15\mu(A\cap T_{c(1)}B\cap C\cap T_{c(1)}^{-1}D),\\
\mu(F^2_{n+1}\cap A\cap T_{l_n}B\cap T_{l_n^2}C\cap T_{l_n^3}D) &\to\frac 15\mu(A\cap B\cap T_{c(1)}C\cap D),\\
\mu(F^3_{n+1}\cap A\cap T_{l_n}B\cap T_{l_n^2}C\cap T_{l_n^3}D) &\to\frac 15\mu(A\cap B\cap C\cap T_{c(1)}D),\\
\mu(F^4_{n+1}\cap A\cap T_{l_n}B\cap T_{l_n^2}C\cap T_{l_n^3}D) &\to\frac 15\mu(A\cap T_{c(1)}^{-1}B\cap T_{c(1)}^{-1}C\cap T_{c(1)}^{-1}D)
\endalign
$$
as $n\to\infty$ with $n\in 3\Bbb N$.
We verify only the first (from the top) claim.
It is straightforward that
$$
\aligned
&\frac{\#(l_nC_{n+1}^0 \triangle c(1)C_{n+1}^1)}{\# C_{n+1}^0}\to 1,\quad \frac{\#(l_nC_{n+1}^1\triangle C_{n+1}^2)}{\# C_{n+1}^1}\to 1,\\
&
\frac{\#(l_nC_{n+1}^2\triangle C_{n+1}^3)}{\# C_{n+1}^2}\to 1,
\quad \frac{\#(l_nC_{n+1}^3\triangle c(1)^{-1}C_{n+1}^4)}{\# C_{n+1}^3}\to 1,\\
&
\frac{\#(l_nC_{n+1}^4\triangle c(1)^{-1}C_{n+1}^0)}{\# C_{n+1}^4}\to 1.
\endaligned
\tag6-1
$$
Take some subsets $\widetilde A,\widetilde B,\widetilde C,\widetilde D\subset F_n\cap c(2)F_n$.
It follows from \thetag{6-1} that
$$
\align
&\mu(F^0_{n+1}\cap [\widetilde A]_n\cap T_{l_n}[\widetilde B]_n\cap T_{l_n^2}[\widetilde C]_n\cap T_{l_n^3}[\widetilde D]_n)\\
=
&\mu([\widetilde AC_{n+1}^0]_{n+1}\cap T_{l_n}[\widetilde BC_{n+1}^4]_{n+1}\cap T_{l_n^2}[\widetilde CC_{n+1}^3]_{n+1}\cap T_{l_n^3}[\widetilde DC_{n+1}^2]_{n+1})+\overline o(1)\\
=
&\mu([\widetilde AC_{n+1}^0]_{n+1}\cap[c(1)^{-1}\widetilde BC_{n+1}^0]_{n+1}\cap[c(1)^{-2}\widetilde CC_{n+1}^0]_{n+1}
\cap[c(1)^{-2}\widetilde DC_{n+1}^0]_{n+1})+\overline o(1)\\
=&\mu([(\widetilde A\cap c(1)^{-1}\widetilde B\cap c(1)^{-2}\widetilde C
\cap c(1)^{-2}\widetilde D)C_{n+1}^0]_{n+1})+\overline o(1)\\
=&\frac 15
\mu([\widetilde A\cap c(1)^{-1}\widetilde B\cap c(1)^{-2}\widetilde C
\cap c(1)^{-2}\widetilde D]_{n})+\overline o(1)\\
=&\frac 15
\mu([\widetilde A]_n\cap T_{c(1)}^{-1}[\widetilde B]_n\cap T_{c(1)}^{-2}[\widetilde C]_n
\cap T_{c(1)}^{-2}[\widetilde D]_{n})+\overline o(1).
\endalign
$$
Hence approximating $A,B,C,D$ with cylinders $[\widetilde A]_n,[\widetilde B]_n,[\widetilde C]_n,[\widetilde D]_n$ respectively and passing to the limit we deduce that
$\mu(F^0_{n+1}\cap A\cap T_{l_n}B\cap T_{l_n^2}C\cap T_{l_n^3}D) \to 0.2\mu(A\cap T_{c(1)}^{-1}B\cap T_{c(1)}^{-2}C\cap T_{c(1)}^{-2}D)$, as claimed.

We now obtain that
$$
\align
\lim_n 5\mu(A\cap T_{l_n}B\cap T_{l_n^2}C& \cap T_{l_n^3}D)
=
\lim_n 5\sum_{j=0}^4\mu(F^j_{n+1}\cap A\cap T_{l_n}B\cap T_{l_n^2}C\cap T_{l_n^3}D)\\
&=\mu(A\cap T_{c(1)}^{-1}B\cap T_{c(1)}^{-2}C\cap T_{c(1)}^{-2}D)\\
&+
\mu(A\cap T_{c(1)}B\cap C\cap T_{c(1)}^{-1}D)
+ \mu(A\cap B\cap T_{c(1)}C\cap D)\\
&+
\mu(A\cap B\cap C\cap T_{c(1)}D)+
\mu(T_{c(1)}A\cap B\cap C\cap D),
\endalign
$$
where the limit is taken along the sequence of $n$ divisible by $3$.
In particular, substituting $A=B=D$ and $C=X$ we obtain
$$
\lim_{n} 5\mu(A\cap T_{l_n}A\cap T_{l_n^3}A) \ge\mu(A)\tag6-2
$$
for each subset $A\subset X$.
On the other hand, take a subset $A\subset X$ such that $\mu(A\cap T_{c(1)}A)=\mu(A\cap T_{c(1)}^2A)=0$.
 Then substituting $B=X$ and $A=C=D$, we obtain
$$
\lim_{n} 5\mu(A\cap T_{l_n}^{-1}A\cap T_{l_n}^{-3}A)=\lim_{n} 5\mu(A\cap T_{l_n^2}A\cap T_{l_n^3}A)=0.\tag6-3
$$
It follows from \thetag{6-2} and \thetag{6-3} that the flows $(T_{c(t)})_{t\in\Bbb R}$
and $(T_{c(-t)})_{t\in\Bbb R}$ are not conjugate.
Thus the second claim of the theorem is proved.

If all $\gamma_n$ are integers, then $l_n$ is power of $c(1)$ and the third claim of the theorem follows from the second one.
\qed
\enddemo

\remark{Remark \rom{6.2}}
In a similar way we can construct  {\it infinite} measure preserving  rank-one actions $T$ of $H_3(\Bbb R)$ for which  the claims of Theorem~6.1 hold.
\endremark

\head 7. Spectral analysis for actions of Heisenberg group
\endhead

Let $T$ be an action of $H_3(\Bbb R)$ on a standard probability space
$(X,\goth B,\mu)$.
Denote by $U$ the corresponding Koopman representation of $H_3(\Bbb R)$
in $L^2(X,\mu)$.
Consider a spectral decomposition of $U$ (we refer to Section~1 for the notation):
$$
\align
 L^2(X,\mu)&= \int^{\oplus}_{\Bbb R^2}\,\bigoplus_{j=1}^{l_U^{1,2}(\alpha,\beta)}\Bbb C\,d\sigma_U^{1,2}(\alpha,\beta)
\oplus\int_{\Bbb R^*}^\oplus\bigoplus_{j=1}^{l_U^3(\gamma)}L^2(\Bbb R,\lambda_\Bbb R)\,d\sigma_U^3(\gamma)\quad\text{and} \\
U(g)&=\int^{\oplus}_{\Bbb R^2}\,\bigoplus_{j=1}^{l_U^{1,2}(\alpha,\beta)} \pi_{\alpha,\beta}(g)\,d\sigma_U^{1,2}(\alpha,\beta)
\oplus\int_{\Bbb R^*}^\oplus\bigoplus_{j=1}^{l_U^3(\gamma)}\pi_\gamma(g)\,d\sigma_U^3(\gamma).
\endalign
$$
We assume that the measures  $\sigma_U^{1,2}$ and $\sigma_U^{3}$ are finite.
Our purpose in this section is to write some easy (almost straightforward) but important corollaries from the spectral decomposition of $U$.

Given $\alpha\in \Bbb R$,
we denote by $\chi_\alpha$ the continuous character $\chi_\alpha(t):=e^{i\alpha t}$ of $\Bbb R$.
We now compute the restriction of $U$ to the subgroup $H_{2,a}$.
For that we identify $H_{2,a}$ with $\Bbb R^2$ via the mapping
$c(t_3)a(t_1)\mapsto (t_3,t_1)$.
It is easy to see that $\pi_{\alpha,
\beta}\restriction H_{2,a}=\chi_0\otimes\chi_\alpha$.
The restriction of $\pi_\gamma$ to $H_{2,a}$ is the tensor product of $\chi_\gamma$ and the left regular representation of $\Bbb R$.
Therefore the `projection' of $\sigma_U^{3}$ to $H_{2,a}=\Bbb R^2$ is equivalent to
the product $\sigma_U^{3}\times
\lambda_{\Bbb R}$.
Hence we obtain the following spectral decomposition for $U\restriction H_{2,a}$:
$$
\align
 L^2(X,\mu)&=
\int^{\oplus}_{\Bbb R}\,\bigoplus_{j=1}^{l_U^{1}(\alpha)}\Bbb C\,d\sigma_U^{1}(\alpha)
\oplus\int_{\Bbb R^*\times\Bbb R}^\oplus\bigoplus_{j=1}^{l_U^3(\gamma)}\Bbb C\,d\sigma_U^3(\gamma)d\lambda_\Bbb R(\alpha)\quad\text{and} \\
U(c(t_3)a(t_1))&=\int^{\oplus}_{\Bbb R}\, \chi_{\alpha}(t_1)I_\alpha\,d\sigma_U^{1}(\alpha)
\oplus\int_{\Bbb R^*\times\Bbb R}^\oplus\chi_\gamma(t_3)\chi_\alpha(t_1)
I_{\gamma}\,d\sigma_U^3(\gamma)d\lambda_\Bbb R(\alpha),
\endalign
$$
where $\sigma^1_U$ is the projection of $\sigma^{1,2}_U$
under the map $\Bbb R^2\ni(\alpha,\beta)\mapsto\alpha\in\Bbb R$, $l_U^1(\alpha)$ is the integral of $l_U^{1,2}$ by the conditional measure $\sigma^{1,2}_U |_{\{\alpha\}\times\Bbb R}$ and
$I_\alpha$ and  $I_\gamma$ are identity operators in
$\bigoplus_{j=1}^{l_U^{1}(\alpha)}\Bbb C$
and $\bigoplus_{j=1}^{l_U^3(\gamma)}\Bbb C$ respectively.

\proclaim{Proposition 7.1}
\roster
\item"\rom(i)"
The  maximal spectral type  of \, $U\restriction H_{2,a}$
 contains the measure $\delta_0\times
\sigma_U^1+\sigma_U^{3}\times
\lambda_{\Bbb R}$ on $\Bbb R^2$.
The corresponding spectral multiplicity map is
$$
\Bbb R^2\ni(\gamma,\alpha)\mapsto
\cases
l^1_U(\alpha) &\text{if $\gamma=0$},\\
l^3_U(\gamma) &\text{otherwise}.
\endcases
$$
\item"\rom(ii)"
The  maximal spectral type  of $U\restriction \{c(t)\mid t\in\Bbb R\}$
 contains the measure $\sigma_U^{1,2}(\Bbb R^2)\delta_0
+\sigma_U^{3}$.
The corresponding spectral multiplicity map is
$$
\Bbb R\ni\gamma\mapsto
\cases
\int l^{1,2}_U\,d\sigma^{1,2} &\text{if $\gamma=0$},\\
\infty &\text{otherwise}.
\endcases
$$
\item"\rom(iii)"
If $(T_{c(t)})_{t\in\Bbb R}$ is ergodic then
$\sigma_U^{1,2}(\Bbb R^2\setminus\{0,0\})=0$, i.e.
there are no non-trivial one-dimensional representations in the spectral decomposition of $U$.
The maximal spectral type of \, $T$ equals the maximal spectral type of
the restriction $T$ to the center of $H_3(\Bbb R)$ (modulo the natural identification).
\item"\rom(iv)"
If $(T_{c(t)})_{t\in\Bbb R}$ is ergodic then $T$ has a simple spectrum if and only if $U\restriction H_{2,a}$  has a simple spectrum.
\item"\rom(v)" If $(T_{c(t)})_{t\in\Bbb R}$ is ergodic then
$(T_{a(t)})_{t\in\Bbb R}$ has Lebesgue (restricted) maximal spectral type.
\item"\rom(vi)"
If $(T_{c(t)})_{t\in\Bbb R}$ is mixing then
$(T_g)_{g\in H_{2,a}}$ is mixing.
\endroster
\endproclaim

We note that (vi) is weaker than  Theorem~3.7.
However we include it here because it follows almost directly from the spectral decomposition for $U$, i.e. it is simpler (shorter) than the proof of Theorem~3.7.

\proclaim{Corollary 7.2} Let $T$ and $S$ be two actions of $H_3(\Bbb R)$
such that the flows $(T_{c(t)})_{t\in\Bbb R}$
and $(S_{c(t)})_{t\in\Bbb R}$
are ergodic.
Then the  set of spectral multiplicities of the Cartesian product $T\times S=(T_g\times S_g)_{g\in H_3(\Bbb R)}$ contains infinity.
\endproclaim
\demo{Idea of the proof}
The claim follows from Propositions~3.10, 7.1(iii) and the following well know fact:
$$
\pi_\gamma\otimes\pi_{\gamma'}=
\cases
\bigoplus_{j=1}^\infty\pi_{\gamma+\gamma'} &\text{if $\gamma+\gamma'\ne 0$},\\
\int_{\Bbb R^2} \pi_{\alpha,\beta}\,d\lambda_\Bbb R(\alpha) \,d\lambda_\Bbb R(\beta) & \text{if $\gamma+\gamma'= 0$}.
\endcases
$$
\qed
\enddemo

It follows, in particular, that the set of spectral multiplicities of each (non-degenerated) Gaussian and Poisson action of $H_3(\Bbb R)$  contains infinity.
The `non-degenerated' here means that that the corresponding (reduced) maximal spectral types are non-atomic.

\head 8. Concluding remarks and open problems
\endhead

\roster
\item
Actions $T$ and $S$ of $H_3(\Bbb R)$ on probability spaces $(X,\mu)$ and $(Y,\nu)$ respectively are called {\it disjoint in sense of Furstenberg \cite{Fu}}
if $\mu\times\nu$ is the only $(T_g\times S_g)_{g\in H_3(\Bbb R)}$-invariant measure on $X\times Y$ with marginals $\mu$ and $\nu$ on $X$ and $Y$ respectively.
Modifying the construction from Section~3 one can obtain  an uncountable family of mutually disjoint
mixing rank-one actions of $H_3(\Bbb R)$ with MSJ.
\item
 The examples of mixing rank-one actions of $H_3(\Bbb R)$ in Section~3 are of stochastic nature.
The choice of the `spacer' maps $s_n$ is not explicit (as in the Ornstein's  example of mixing $\Bbb Z$-action \cite{Or}).
Is it possible to construct  explicit, concrete examples of such actions
by analogy with the classical $\Bbb Z$-staircases \cite{Ad}?
\item
  Are there smooth (differentiable) models for mixing rank-one actions of
$H_3(\Bbb R)$?
\item
All the results of this paper obtained for the 3-dimensional Heisenberg group extend naturally to the Hisenberg groups $H_{2k+1}(\Bbb R)$, $k>1$, of higher dimensions and some `generalized' Heisenberg groups \cite{Ki}.
Is it possible to extend them to   the class of all (or a sufficiently wide subclass of) simply
connected nilpotent Lie groups?
\item
We conjecture that a mixing rank-one action of $H_3(\Bbb R)$ is mixing of all orders.
\item
 Whether each 2-fold simple weakly mixing action of $H_3(\Bbb R)$ is simple?
\item
Let an action $T$ of $H_3(\Bbb R)$ have MSJ.
Does this imply that the flow $(T_{c(t)})_{t\in\Bbb R}$ is simple?
Is this true in the particular case when $T$ is of rank one?
\item
Can we construct an action $T$ of $H_3(\Bbb R)$ such that
the flow $(T_{c(t)})_{t\in\Bbb R}$ is ergodic and conjugate to its inverse but the flows $(T_{a(t)})_{t\in\Bbb R}$ and
$(T_{b(t)})_{t\in\Bbb R}$
are non-conjugate or even disjoint in the sense of Furstenberg?
We note that  $T$ is always `spectrally' symmetric, i.e. the Koopman representation $U_T$ is unitarily equivalent to $U_T\circ\theta$, where $\theta$ denotes the flip in $H_3(\Bbb R)$.
Indeed, $\pi_\gamma\circ\theta$ is unitarily equivalent to $\pi_{-\gamma}$ for each $\gamma\in\Bbb R^*$ and the measure $\sigma_U^3$ is quasi-invariant under the inversion $\Bbb R^*\ni t\mapsto -t\in\Bbb R^*$ in view of Proposition~7.1(iii). 
In particular, we  now deduce from  Proposition~7.2(v) 
 $(T_{a(t)})_{t\in\Bbb R}$ and
$(T_{b(t)})_{t\in\Bbb R}$
have the same Lebesgue maximal spectral type.
\item
Is there an action $T$ of $H_3(\Bbb R)$ such that
the flow $(T_{c(t)})_{t\in\Bbb R}$ is ergodic and conjugate to its inverse,
the flows $(T_{a(t)})_{t\in\Bbb R}$ and
$(T_{b(t)})_{t\in\Bbb R}$
are conjugate but $T$ is asymmetric, i.e. it is not conjugate to the action 
$(T_{\theta(g)})_{g\in H_3(\Bbb R)}$?
\item
Suppose that $T$ and $S$ are two disjoint ergodic actions  of $H_3(\Bbb R)$.
Does this imply that the flows $(T_{c(t)})_{t\in\Bbb R}$
and
$(S_{c(t)})_{t\in\Bbb R}$ are also disjoint?
 \endroster

\enddocument